\documentclass[USenglish, a4paper, oneside, 10pt]{article}
\makeatletter
\let\stdmaketitle\maketitle
\makeatother

\PassOptionsToPackage{capitalize}{cleveref}
\PassOptionsToPackage{noend}{algpseudocode}

\def\embeddedpreamble{myPreamble}
\usepackage[
	arxiv,
	alglinenumber,
	eqreset=section,
	thmreset=section,
]{\embeddedpreamble}
\providecommand{\keywords}[1]{\textbf{Keywords:} #1}
\providecommand{\amsclassification}[1]{\textbf{2020 AMS Subject Classification:} #1}
\usepackage{amsmath}
\usepackage{amssymb}
\usepackage{amsfonts}
\usepackage[labelsep=colon]{caption}
\usepackage{stackrel}

\usepackage{url}
\usepackage{multirow}

\usepackage{subcaption}
\usepackage{pgf}
\usepackage{pgfplots}
\usepgfplotslibrary{fillbetween}
\usepackage{tikz}
\usetikzlibrary{arrows,automata}
\usetikzlibrary{intersections}
\usetikzlibrary{calc}
\usetikzlibrary{shapes}
\usetikzlibrary{positioning}
\usetikzlibrary{decorations.text}
\pgfdeclarelayer{background}
\pgfdeclarelayer{foreground}
\pgfsetlayers{background,main,foreground}
\usetikzlibrary{calc,patterns,decorations.pathmorphing,decorations.markings}

\usepackage{xcolor}

\definecolor{mycol}{RGB}{19,48,128}

\definecolor{reddishbrown}{HTML}{A52A2A}

\definecolor{PaperRed}{RGB}{204,0,0}
\definecolor{PaperGray}{gray}{0.55}

\definecolor{PaperBlue}{RGB}{0,90,180}
\definecolor{PaperGreen}{RGB}{0,120,60}
\newif\ifshownew    \shownewtrue
\newif\ifshowold    \showoldtrue
\newif\ifshowrev    \showrevtrue
\newif\ifdiffcolor  \diffcolortrue

\graphicspath{{./}}

\newcommand{\NN}{{\mathcal N}}

\newcommand{\rr}{\R}

\newcommand{\beqar}{\begin{eqnarray}}
	\newcommand{\eeqar}{\end{eqnarray}}
\newcommand{\beqarno}{\begin{eqnarray*}}
	\newcommand{\eeqarno}{\end{eqnarray*}}
\newcommand{\ba}[1]{\begin{array}{#1}}
	\newcommand{\ea}{\end{array}}
\newcommand{\matrice}[2]{\left[\hspace*{-.1cm}\ba{#1} #2 \ea\hspace*{-.1cm}\right]}

\newcommand{\algnamefont}{\ttfamily}
\newcommand{\BALM}{{\hyperref[alg:BALM]{\bf \algnamefont FS-ALM}}}
\newcommand{\refBALM}{\BALM~(\cref{alg:BALM})}

\newcommand{\PBALM}{{\hyperref[alg:PBALM]{\bf \algnamefont FS-P-ALM}}}

\usepackage{graphicx}
\usepackage{amsthm}
\usepackage{algorithm}
\usepackage{commath}
\usepackage{color}
\usepackage{mathtools}
\usepackage{graphicx}
\usepackage{booktabs}
\usepackage{adjustbox}
\usepackage{tabularx}
\usepackage{float}
\usepackage{titlesec}

\usepackage{enumitem}

\newcommand{\Secref}[1]{Section \ref{#1}}

\allowdisplaybreaks

\newcommand{\calK}{{\mathcal{K}}}

\newcommand{\calN}{{\mathcal{N}}}
\newcommand{\calS}{{\mathcal{S}}}

\newcommand{\calL}{{\mathcal{L}}}

\newcommand{\rnn}{{\rr_{+}}}
\newcommand{\rp}{{\rr_{++}}}

\newcommand{\gradn}{\mathop{\nabla}\nolimits}

\usepackage{accents}

\makeatletter

	\newcommand{\Bk}{\@ifstar{B_{k+1}}{B_k}}
	\newcommand{\Ck}{\@ifstar{\C_{k+1}}{\C_k}}
	\newcommand{\fk}{\@ifstar{\tilde f_{k+1}}{\tilde f_k}}
	\newcommand{\gk}{\@ifstar{\tilde g_{k+1}}{\tilde g_k}}
	\newcommand{\hk}{\@ifstar{{\red h_{k+1}}}{{\red h_k}}}
	\newcommand{\lk}{\@ifstar{{\red\ell_{k+1}}}{{\red\ell_k}}}
	\newcommand{\Lk}{\@ifstar{{\red L_{k+1}}}{{\red L_k}}}
	\newcommand{\LLk}{\@ifstar{{\red L_{k+1}}}{{\red L_k}}}
	\newcommand{\lamk}{\@ifstar{\lambda_{k+1}}{\lambda_k}}
	\newcommand{\Lamk}{\@ifstar{\Lambda_{k+1}}{\Lambda_k}}
	\newcommand{\Pk}{\@ifstar{P_{k+1}}{P_k}}
	\newcommand{\Uk}{\@ifstar{\U_{k+1}}{\U_k}}

	\newcommand{\alphk}{\@ifstar{\alpha_{k+1}}{\alpha_k}}
	\newcommand{\betk}{\@ifstar{\beta_{k+1}}{\beta_k}}
	\newcommand{\Delk}{\@ifstar{\Delta_{k+1}}{\Delta_k}}
	\newcommand{\epsk}{\@ifstar{\varepsilon_{k+1}}{\varepsilon_k}}
	\newcommand{\gamk}{\@ifstar{{\gamma_{k+1}}}{{\gamma_k}}}

	\newcommand{\bgamk}{\@ifstar{{\bar\gamma}_{k+1}}{{\bar\gamma}_k}}

	\newcommand{\Hk}{\@ifstar{H_{k+1}}{H_k}}
	\newcommand{\rhok}{\@ifstar{\rho_{k+1}}{\rho_k}}
	\newcommand{\nuk}{\@ifstar{\nu_{k+1}}{\nu_k}}
	\newcommand{\sigk}{\@ifstar{\sigma_{k+1}}{\sigma_k}}
	\newcommand{\xik}{\@ifstar{\xi_{k+1}}{\xi_k}}

	\newcommand{\Fwk}{\@ifstar\@@Fwk\@Fwk}
	\newcommand{\Fw@twoargs}[2][]{{\nocolor\ifstrempty{#2}{\id-\gamma_{#1}\nabla f}{#2-\gamma_{#1}\nabla f(#2)}}}
	\newcommand{\@Fwk}[2][k]{\Fw@twoargs[#1]{#2}}
	\newcommand{\@@Fwk}[2][k+1]{\Fw@twoargs[#1]{#2}}

	\newcommand{\FBk}{\@ifstar\@@FBk\@FBk}
	\newcommand{\FB@twoargs}[2][]{{\nocolor\prox_{\gamma_{#1}g}(\Fw@twoargs[#1]{#2})}}
	\newcommand{\@FBk}[2][k]{\FB@twoargs[#1]{#2}}
	\newcommand{\@@FBk}[2][k+1]{\FB@twoargs[#1]{#2}}

	\newcommand{\PDx}{\@ifstar\@@PDx\@PDx}
	\newcommand{\PDx@threeargs}[3][]{{\nocolor\prox_{\gamma_{#1}g}\left(\Fw@twoargs[#1]{#2}-\gamma_{#1}\linop*\ifstrempty{#3}{{}\cdot{}}{#3}\right)}}
	\newcommand{\@PDx}[3][k]{\PDx@threeargs[#1]{#2}{#3}}
	\newcommand{\@@PDx}[3][k+1]{\PDx@threeargs[#1]{#2}{#3}}

	\def\innprod{\@ifstar\@innprod\@@innprod}
	\newcommand{\linop}{\@ifstar{\trans A}{A}}

	\newcommand{\C}{c}

	\newcommand{\U}{\mathcal U}

	\newcommand{\DRE}{\varphi^{\text{\sc dr}}\@ifnextchar_{}{_\gamma}}

	\newcommand{\disablecolorlinks}{\def\HyColor@UseColor##1{}}

\showoldfalse
\diffcolorfalse

\makeatother

\makeatletter
\let\maketitle\stdmaketitle
\makeatother

\title{
	A proximal augmented Lagrangian method for nonconvex optimization with equality and inequality constraints
	\thanks{
		This work was funded by the European Union (ERC Advanced Research Grant COMPACT, No. 101141351). Views and opinions expressed are however those of the authors only and do not necessarily reflect those of the European Union or the European Research Council. Neither the European Union nor the granting authority can be held responsible for them.
		P. Latafat is a member of the Gruppo Nazionale per l'Analisi Matematica, la Probabilit\`a e le loro Applicazioni (GNAMPA - National Group for Mathematical Analysis, Probability and their Applications) of the Istituto Nazionale di Alta Matematica (INdAM - National Institute of Higher Mathematics)
	}
}

\author{
	Adeyemi D.~Adeoye\thanks{
		DYSCO (Dynamical Systems, Control, and Optimization), IMT School for Advanced Studies Lucca, 
		Piazza S.Francesco, 19, 55100 Lucca, Italy.
		{\it E-mails:} \{adeyemi.adeoye,~puya.latafat,~alberto.bemporad\}@imtlucca.it}
	\and
	Puya Latafat\footnotemark[2]
	\and
	Alberto Bemporad\footnotemark[2]
}
\date{}

\begin{document}
\maketitle

\begin{abstract}
	We propose an inexact proximal augmented Lagrangian method (P-ALM) for solving nonconvex structured optimization problems. The proposed method features an easily implementable rule not only for updating the penalty parameters, but also for adaptively tuning the proximal term. It allows the penalty parameter to grow rapidly in the early stages to speed up progress, while ameliorating the issue of ill-conditioning in later iterations, a well-known drawback of the traditional approach of linearly increasing the penalty parameters. A key element in our analysis lies in the observation that the augmented Lagrangian can be controlled effectively along the iterates, provided an initial feasible point is available. Our analysis, while simple, provides a new theoretical perspective about P-ALM and, as a by-product, results in similar convergence properties for its non-proximal variant, the classical augmented Lagrangian method (ALM). Numerical experiments, including convex and nonconvex problem instances, demonstrate the effectiveness of our approach.
\end{abstract}
\\[1ex]
\amsclassification{Primary
65K05,
93-08,
49M37;
Secondary
90C06,
90C53.
}

\section{Introduction}\label{sec:intro}
We consider the constrained nonlinear programming problem
\begin{equation}
    \label{eq:prob}
    \begin{aligned}
        \minimize_{x\in\rr^n}&\quad f(x) \coloneqq f_1(x) + f_2(x) \\
        \stt{}& \quad h(x) = 0 \\
        &\quad g(x) \leq 0,
    \end{aligned}
\end{equation}
where $\func{f_1}{\R^n}{\rr}$, $\func{h}{\R^n}{\R^p}$ and $\func{g}{\R^n}{\R^m}$ are continuously differentiable and possibly nonconvex functions, and $\func{f_2}{\R^n}{\Rinf}$ is a proper closed convex function whose proximal mapping is easy to evaluate. Problem \eqref{eq:prob} covers many optimization problems in engineering, finance, and machine learning. Typically, $f_1$ represents an objective function, $h$ and $g$ equality and inequality constraints on decision variables and $f_2$ can capture regularization terms, such as the $1$-norm, or encodes simple constraints, like box constraints. 
Because of this flexible structure, formulations of this type routinely arise across diverse domains, including optimal control \cite{mayne2000constrained}, statistical machine learning \cite{hastie2009elements}, and operations research \cite{boyd2004convex}.

The augmented Lagrangian method (ALM), also known as the method of multipliers,  first introduced by Hestenes \cite{hestenes1969multiplier} and Powell \cite{powell1969method}, is a widely used approach for solving constrained optimization problems of the form \eqref{eq:prob}. In \cite{rockafellar1976augmented}, Rockafellar proposed an elegent convergence analysis for ALM by viewing it as an instance of the proximal point algorithm (PPA) \cite{rockafellar1976monotone} applied to an operator derived from the dual problem. Specifically, standard ALM updates take the form 
\cite{rockafellar1976monotone,rockafellar1976augmented}:
    \begin{align}
        x^{k+1} \in{}&
        \argmin_{x\in\rr^n} \left\{f(x) + \langle \lambda^k, h(x) \rangle + \tfrac{\rho}{2} \|h(x)\|^2 + \tfrac{1}{2\nuk} \|\max\{0, \nuk g(x) + \mu^k\}\|^2 
        \right\}
        \label{eq:x-update}\\
        \lambda^{k+1} &= \lambda^k + \rho_k h(x^{k+1})\notag\\
        \mu^{k+1} &= \max\{0, \mu^k + \nu_k g(x^{k+1})\},\notag
    \end{align}
where \(\rho,\nu\in\rp\) are penalty parameters, and $\lambda\in\rr^p$, $\mu\in\rr_+^m$ are the Lagrange multipliers for equality and inequality constraints, respectively.

The primary focus of this paper is on the proximal variant of ALM, known as the proximal augmented Lagrangian method (P-ALM). Rockafellar proposed P-ALM as an instance of PPA applied to the operator related to a primal-dual formulation, see \eqref{eq:proxALM-fixed-point}. The iterates of proximal ALM differs from the ALM formulation in the update of the primal variable. In particular, \eqref{eq:x-update} is replaced with 
    \begin{align}\label{eq:proxALM-updates}
        x^{k+1} &\in \argmin_{x\in\rr^n} \left\{\calL_{\rho_k,\nu_k,\gamma_k}(x, \lambda^k,\mu^k; x^k) + f_2(x)\right\},
    \end{align}
where
\begin{align}
    \calL_{\rho,\nu,\gamma}(x, \lambda, \mu; v) 
    \coloneqq{}& f_1(x) + \langle \lambda, h(x) \rangle + \frac{\rho}{2} \|h(x)\|^2 {}+ \frac{1}{2\nu} \|\max\left\{0, \nu g(x) + \mu\right\}\|^2 \nonumber 
    \\
    &
    {} - \frac1{2\nu}\|\mu\|^2
    {}+ \frac{1}{2\gamma}\|x - v\|^2,\label{eq:prox-ALM:lagrangian}
\end{align}
is the proximal augmented Lagrangian function associated with problem \eqref{eq:prob}.
Note that we intentionally isolate $f_2$ from the definition of proximal augmented Lagrangian as it will be handled without introducing additional multipliers, provided its proximal mapping can be evaluated efficiently. Two distinct uses of ``proximal'' are at play here. The regularization term $\frac{1}{2\gamma}\|x-v\|^2$ in \eqref{eq:prox-ALM:lagrangian} is what distinguishes P-ALM from ALM and is what the prefix ``P-'' refers to throughout, whereas the proximal mapping of $f_2$ just mentioned is how the nonsmooth term is handled by the inner solver and is present in both variants; setting $\gamma_k\equiv\infty$ removes the former and recovers ALM, leaving the latter unaffected.

\subsection{Motivation and contributions}

Recent years have witnessed a renewed interest in P-ALM and its non-proximal variant, ALM. Initially popular because they could efficiently solve a wide range of optimization problems, they were eclipsed by interior-point methods; however, their superior scalability on very large-scale, highly structured problems has spurred a resurgence among researchers seeking alternatives to conventional interior-point and sequential quadratic programming methods \cite{shefi2014rate,curtis2015adaptive}. In this work, we focus our convergence analysis primarily on P-ALM, with the results for ALM derived as a by-product.

A major advantage of the P-ALM, over its non-proximal variant, is that it results in Lagrangian subproblems that are better behaved in general. For instance, depending on the proximal parameter, a nonconvex subproblem may become (strongly) convex, or a lower unbounded subproblem may be guaranteed to admit a solution.
Moreover, the proximal variant is closely related to many popular numerical methods such as the alternating direction method of multipliers \cite{douglas1956numerical, glowinski1989augmented, tseng1991applications, eckstein1992douglas}, predictor-corrector proximal multiplier \cite{chen1994proximal}, and primal-dual hybrid gradient method \cite{chambolle2011first}. We refer the reader to  \cite{bertsekas2014constrained,shefi2014rate} for a complete account of developments and variants of augmented Lagrangian methods. 

The classical convergence analysis for augmented Lagrangian type methods crucially hinges on the maximal monotonicity assumption for the associated operator \cite{rockafellar11970monotone,eckstein1992douglas}. Such analysis readily accommodates inexact variants under summability of the error sequence. The inexact case is particularly relevant for large-scale or structured problems, where solving subproblems exactly is impractical. In the monotone setting, relatively inexact variants of ALM have also been studied eliminating the need for a predefined sequence of summable inexactness tolerance while ensuring convergence (see \cite{eckstein2013practical}).

However, these powerful ideas do not extend to the nonconvex setting studied here. Instead, we adopt an inexactness framework as in \cite{birgin2018augmented,birgin2020complexity}, which involves finding an approximate stationary point of the (proximal) augmented Lagrangian subproblem up to a fixed tolerance. Consequently, our analysis establishes termination at an approximate KKT point (see \cref{def:epsKKT}) rather than proving sequential convergence.

\begin{figure}[H]
	\centering
	\subfloat{
		\resizebox*{6cm}{!}{\includegraphics{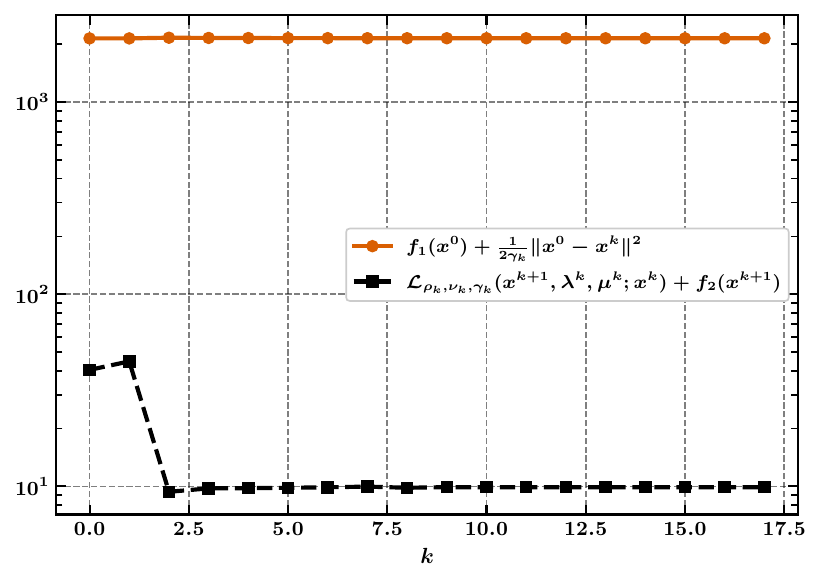}}}
	\hspace{0.5cm}
	\subfloat{
		\resizebox*{6cm}{!}{\includegraphics{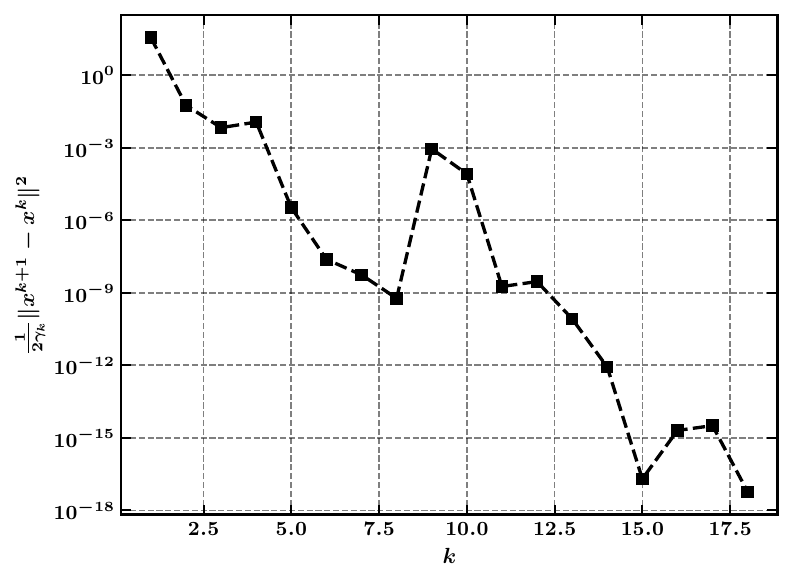}}}
	\caption{Illustration of two defining properties of \PBALM{} central to the proof of \cref{thm:ALM:convergence}, on the basis-pursuit instance of \Secref{sec:basis_pursuit} ($p=400$, $n=1024$), with $\alpha=1.01$ under \texttt{PANOC}. \textbf{Left:} the augmented Lagrangian value at $x^{k+1}$ stays below the reference bound at every iteration. \textbf{Right:} the proximal term vanishes.}
	\label{fig:lagrangian-bound}
\end{figure}

Our convergence analysis is largely inspired by recent results for ALM  in \cite{grapiglia2021complexity,evens2021neural}. A common feature of these works, shared also with ours, is the assumption that an initial feasible point is known, which is used to ensure the boundedness of the augmented Lagrangian function (see \cite[Ex. 4.12]{birgin2014practical}; see also \cref{fig:lagrangian-bound} for a practical illustration). While \cite{grapiglia2021complexity} considers both equality and inequality constraints, it implicitly requires the Lagrange multipliers to converge to zero for the penalty parameters to remain bounded. This limitation arises from the specific way the penalty parameters are updated, involving a scaled norm of the multiplier vectors \cite[Eqs. (2.8) and (2.9)]{grapiglia2021complexity}. This issue was resolved in \cite{evens2021neural} for problems with equality constraints within an ALM framework.
This mechanism is different from the multiplier safeguarding used in \cite{de2023constrained,jia2023augmented} and, classically, in \cite{birgin2014practical}. Rather than projecting or clipping the multipliers to a prescribed bounded set, we control their growth as a by-product of the initial feasibility condition. In fact,  $\lambda^k,\mu^k$ need not stay bounded, but their growth remains tied to that of the penalty parameters (\cref{thm:proxALM:multipliers_bound}), which is what the finite-time termination argument uses.
Building on these works, we provide a new analysis for P-ALM addressing \eqref{eq:prob} in its full generality, with ALM recovered as a special case. Unlike \cite{evens2021neural}, which treats equality constraints only, we also handle inequalities. The analysis is carried out for the proximal variant directly, with ALM recovered as the special case $\gamma_k\equiv\infty$, so the growth conditions on $\gamma_k$ in \cref{assump:penalty} must interact correctly with the penalty schedule. Moreover, the nonsmooth term $f_2$ enters our notion of $\varepsilon$-KKT point (\cref{def:epsKKT}) through a proximal residual rather than a gradient, so the stationarity measure the inner solver must certify is correspondingly different.

\subsubsection{Adaptive parameter updates}\label{sec:adaptive-parameter}
A crucial practical issue in augmented Lagrangian methods lies in the update of the penalty parameters \(\rhok, \nuk\). The standard approach is to increase the penalty parameters linearly, \eg, \(\rhok* = \xi \rhok\), where $\xi > 1$, whenever the quantity \(\|h(x^{k+1})\|_{\infty}\) is not improved compared to the prior iterate (and similarly for the penalty parameter \(\nuk\) associated with the inequality contraints). 
However, it is well-known that a linear rate of increase for the penalty parameter can lead to ill-conditioning \cite{bertsekas1975combined, bertsekas1976multiplier}. Recent nonconvex augmented Lagrangian methods with a global convergence guarantee \cite{chen2017augmented,de2023constrained,jia2023augmented}, as well as the classical safeguarded framework of \cite{birgin2014practical}, all update the penalty this same multiplicative way. In contrast, whenever the feasibility measure does not improve, our analysis allows for parameter updates of the form
\begin{align}
    \rhok* = \max\set{\rhok, \rho_0 (k+1)^\alpha},\notag
\end{align}
where the penalty parameter grows rapidly in the initial phases but asymptotically grows at a much lower rate than linear growth, thereby ameliorating the issue of ill-conditioning. The penalty parameter \(\nuk\) is updated similarly (see our \PBALM{} algorithm below for details). Throughout, the prefix \textsc{fs} stands for \emph{feasible start} as the methods assume that a feasible $x^0$ is available. This is crucial for the control of the augmented Lagrangian and of the multiplier sequence (see \cref{thm:proxALM:L-bound,thm:eta-vanish}). 

Moreover, our analysis suggests an update rule for the proximal parameter $\gamma_k$ presented in \cref{state:proxALM:eta} of \PBALM{} (see also \cref{assump:penalty}). This results in stronger regularization in early iterations while gradually transitioning the method toward a standard ALM variant.

\subsection{Other related work}

While our primary focus is the P-ALM, we highlight that several notable works have studied the classical ALM in the convex setting beyond the traditional PPA framework, particularly with respect to iteration-complexity estimates, see \cite{lan2016iteration, patrascu2017adaptive, necoara2019complexity, birgin2020complexity, li2021inexact, xu2021iteration, lu2023iteration}. For nonconvex problems with equality constraints \cite{sahin2019inexact, li2021rate, bodard2024inexact} employed a regularity condition to study the global convergence and worst-case iteration-complexity bounds for ALM. The recent work \cite{liu2025lower} derived lower complexity bounds for certain classes of first-order methods that solve linearly-constrained nonconvex problems, including ALM. 
ALM variants in which simple constraints such as boxes or polytopes are handled separately from the equality and inequality constraints were studied in \cite{conn1991globally, conn1996convergence, andreani2008augmented, birgin2020complexity}. 
For P-ALM, this formulation was studied in \cite{zhang2020proximal, zhang2022global}, where iteration-complexity bounds were derived for nonconvex problems with linear constraints.

Despite its extensive practical use, P-ALM has not yet been fully understood for nonconvex problems. Recent studies such as \cite{zhang2022global, xie2021complexity}, derive outer-complexity estimates for special cases of \eqref{eq:prob}, such as linearly-constrained and/or equality-constrained problems. Because of its close connection to the alternating direction method of multipliers, other seemingly-tangential but important results were developed in \cite{hong2017prox, jiang2019structured, zhang2020proximal} for problems with linear equality constraints, and in \cite{hajinezhad2019perturbed, boct2020proximal} where the objective may be nonsmooth and the penalty parameter is kept constant. The authors in \cite{sujanani2023adaptive, kong2023iteration} focused on variants of the P-ALM where an accelerated gradient method solves nonconvex subproblems and demonstrates improved computational performance for linearly-constrained problems. Finally, \cite{hermans2022qpalm} developed an inexact P-ALM for quadratic programming (QP) problems with linear constraints.

In contrast to most of these works, we present a general P-ALM framework that handles nonlinear equality and inequality constraints, as well as nonsmooth terms. Despite its simplicity, our analysis covers both the more stabilized Proximal variant (P-ALM) and the standard ALM settings while suggesting practical rules for updating the penalty parameters and adaptively tuning the proximal term. Additionally, our approach guarantees explicit control of the augmented Lagrangian and the multipliers under an initial feasibility condition (see \cref{thm:proxALM:L-bound}). This feature is central to our main convergence results presented in \cref{thm:ALM:convergence}. Overall, our results provide new theoretical guarantees and practical guidelines for designing efficient P-ALM algorithms applicable to a wide class of constrained optimization problems.

\subsection{Organization}
The rest of this paper is organized as follows. In \Secref{sec:notations}, we introduce the notation conventions used throughout the paper. Then, in \Secref{sec:preliminaries}, we introduce the main definitions and problem assumptions adopted throughout the paper. \Secref{sec:algorithms} describes the proposed P-ALM (\PBALM) and presents its main convergence properties. In \Secref{sec:convergence}, we establish several key theoretical properties of \PBALM{} and present the proof of the main convergence results. We present an ALM (\BALM) as a special case of \PBALM{} in \Secref{sec:ALM-algorithm} and establish its key theoretical properties from those of \PBALM. Finally, \Secref{sec:experiments} reports results from numerical simulations to illustrate the performance of our methods on both convex and nonconvex problems, including instances from the Maros-M{\'e}sz{\'a}ros collection \cite{maros1999repository} and the basis pursuit problem \cite{chen2001atomic} on synthetic datasets.

\subsection{Notations}\label{sec:notations}
The set of real and extended-real numbers are \(\R\coloneqq(-\infty,\infty)\) and \(\Rinf\coloneqq\R\cup\set\infty\), while the positive and strictly positive reals are \(\R_+\coloneqq[0,\infty)\) and \(\R_{++}\coloneqq(0,\infty)\).
We use the notation \([x]_+ = \max\set{0, x}\).
With \(\id\) we indicate the identity function \(x\mapsto x\) defined on a suitable space.
We denote by $\R^n$ the standard $n$-dimensional Euclidean space with inner product \(\innprod{{}\cdot{}}{{}\cdot{}}\) and induced norm \(\|{}\cdot{}\|\). For a fixed integer \(n\ge 1\) we write  $\R_+^n \coloneqq [0,\infty)^n=\{x=(x_{1},\dots ,x_{n})\in\R^n \mid x_i\in [0,\infty) \text{ for all } i=1,\ldots ,n\}$.
For a vector $w = (w_1, \hdots, w_N) \in \R^n$, $w_i \in \R^{n_i}$ is used to denote its $i$-th (block) coordinate.
The infinity norm, and $1$-norm are denoted by $\|\cdot\|_\infty$, and $\|\cdot\|_1$, respectively.
We use the notation 
$
    \seq{
        z^k
        }[k\in I]
$
to denote a sequence with indices in the set $I\subseteq \N$. When dealing with scalar sequences we use the subscript notation $\seq{\gamma_k}[k\in I]$. The set of indices $\{1, 2, \ldots, m\}$ is denoted by $[m]$.

An operator or set-valued mapping $A:\R^n\rightrightarrows\R^d$ maps each point $x\in\R^n$ to a subset $A(x)$ of $\R^d$. 
We denote the domain of $A$ by $\dom A\coloneqq\{x\in\R^n\mid A(x)\neq\emptyset\}$. For an extended real-valued function $f : \R^n \rightarrow \Rinf$, the domain is defined as $\dom f \coloneqq \set{x \in \R^n}[f(x) < \infty]$. We say that $f$ is proper if $\dom f \neq \emptyset$ and that $f$ is closed or lower semicontinuous if its epigraph $\epi f \coloneqq \set{(x,\alpha) \in \R^n \times \R}[f(x) \leq \alpha]$ is a closed subset of $\R^{n+1}$. 
For a differentiable function \(f\), $\gradn_{x} f$ denotes its gradient with respect to \(x\). For a convex function \(f\), $\partial_xf $ is used to denote the subdifferential of $f$ with respect to $x$. When the differentiation is taken with respect to the full variable (i.e., all components), we omit the subscript and simply write $\partial f$. For a vector-valued function $f\colon\R^n\to\R^m$, we denote the Jacobian of $f$ at $x$ by $J_f(x)$.
The indicator function of a set $C\subseteq\R^n$ is denoted by \(\indicator_C\), namely \(\indicator_C(x)=0\) if \(x\in C\) and \(\indicator_C(x) =\infty\) if \(x \notin C\). We denote the normal cone of $C$ by \(N_C\). 
The proximal operator of a function $f$ at $x$ is defined as
\begin{equation*}
    \prox_{\gamma f}(x) \coloneqq \argmin_{y\in\R^n} \left\{f(y) + \tfrac{1}{2\gamma}\|y-x\|^2\right\}, \quad \forall x\in\R^n, \gamma \in \rp.
\end{equation*}
\section{Assumptions and preliminaries}\label{sec:preliminaries}
In the following, we present the key assumptions we use throughout the paper.
\begin{assumption}[Problem assumptions]\label{assump:prob}
    The following conditions hold for problem \eqref{eq:prob}:
    \begin{enumeratass}
    	\item The functions $\func{f_1}{\R^n}{\rr}$, $\func{h}{\R^n}{\R^p}$, and $\func{g}{\R^n}{\R^m}$ are continuously differentiable.
        \item The function $\func{f_2}{\R^n}{\Rinf}$ is proper, closed and convex, and its proximal mapping is easy to compute. \label{assump:f2}
        \item There exists $\bar{x} \in \dom f_2$ such that $h(\bar{x}) = 0$ and $g(\bar{x}) \leq 0$. \label{assump:feas_sol}
        \item The objective function $f$ is lower-bounded by a finite number, \ie, $-\infty < f^\star  \coloneqq \inf f$ for some $f^\star \in \R$. \label{assump:obj}
    \end{enumeratass}
\end{assumption}
The lower boundedness of the objective function required in \cref{assump:obj} ensures the well-definedness of the Lagrangian subproblems, as it guarantees the existence of \(\varepsilon\)-stationary points required in \eqref{prop:proxALM:grad}.

A tuple $(x,\lambda,\mu)$ of primal and dual variables satisfies the Karush-Kuhn-Tucker (KKT) conditions for problem \eqref{eq:prob} if:
\begin{equation}
    \label{eq:KKT}
    \begin{aligned}
    & 
        \|x - \prox_{f_2}(x - \nabla_{x}\calL(x,\lambda,\mu))\| = 0, \\
    &
        h(x) = 0, \quad 
        g(x) \leq 0, \\
    & 
        \mu_i \geq 0, \quad 
        \mu_i g_i(x) = 0 \quad \text{for all } i \in [m],
\end{aligned}
\end{equation}
where the Lagrangian function $\calL$ is defined by
\begin{align}
    \calL(x,\lambda,\mu) \coloneqq f_1(x) + \langle\lambda, h(x)\rangle + \langle\mu, g(x)\rangle.
    \label{eq:L}
\end{align} 
Here, the primal stationarity condition is expressed through the natural residual mapping \(x - \prox_{f_2}(x - \nabla_{x}\calL(x,\lambda,\mu))\) that vanishes if and only if \(-\nabla_x \calL(x, \lambda, \mu) \in \partial f_2(x)\).

Finding a point that satisfies the KKT conditions is equivalent to finding a zero of the primal-dual (KKT) operator $T\colon\R^n \times\R^p\times\R^m\rightrightarrows\R^n\times\R^p\times\R^m$ given by
\begin{align*}
    T(x,\lambda,\mu) \coloneqq \matrice{c}{
        \nabla f_1(x) + \partial f_2(x) + J_h(x)^\top \lambda + J_g(x)^\top\mu\\
        -h(x) \\
        -g(x) + \NN_{\R_+^m}(\mu)}.
\end{align*}
That is, the first order optimality conditions of P-ALM updates in \eqref{eq:proxALM-updates} can be written compactly as
\begin{equation}\label{eq:proxALM-fixed-point}
    \begin{aligned}
        (\tfrac{1}{\gamma_k}(x^k - x^{k+1}), \tfrac{1}{\rho_k}(\lambda^k-\lambda^{k+1}), \tfrac{1}{\nu_k}(\mu^k-\mu^{k+1})) \in T(x^{k+1},\lambda^{k+1},\mu^{k+1}).
    \end{aligned}
\end{equation}
In practice, one aims to solve these fixed-point equations approximately. To this end, we consider the following notion of $\varepsilon$-KKT optimality for P-ALM subproblems.

\begin{definition}[$\varepsilon$-KKT optimality]
    \label{def:epsKKT}
    A point $(x,\lambda,\mu)$ is called an $\varepsilon$-KKT point of \eqref{eq:prob} if it satisfies the following conditions:
    \begin{equation*}
        \begin{aligned}
            &{}\|x - \prox_{f_2}(x - \nabla_{x}\calL(x,\lambda,\mu))\|_\infty
        {}\leq{}
            \varepsilon,\quad 
        {}
            \|h(x)\|_\infty
        {}\leq{}
            \varepsilon, \quad
            \|\max\{0,g(x)\}\|_\infty \leq \varepsilon,
            \\
            &{}
            \mu_i \geq 0, \quad    
            \mu_i = 0\; 
            \text{ whenever } g_i(x) < -\varepsilon, \text{ for all } i \in [m] .
        \end{aligned}
    \end{equation*}
\end{definition}
The $\varepsilon$-KKT conditions are a relaxation of the standard KKT conditions, allowing for approximate satisfaction of the optimality conditions. They are akin to those of \cite{birgin2020complexity} which studies the specific case where $f_2$ is the indicator function of a box.

\section{An inexact proximal ALM}\label{sec:algorithms}

Consider the inexact proximal augmented Lagrangian method \PBALM{} (defined in \cref{alg:PBALM}). Given an initial feasible point $x^0$, the algorithm sets the penalty parameters $\rho_0$ and $\nu_0$, and the associated Lagrange multipliers which are iteratively updated in \cref{state:ALM:lambda,state:ALM:mu,state:proxALM:rho,state:proxALM:nu}. As is common in the literature, \PBALM{} probes the changes in the quantity \(\|h(x^k)\|_{\infty}\) along iterations to decide if the penalty parameter corresponding to the equality constraints should be increased or not. For inequality constraints, it employs the quantity $E^k$ defined by 
\begin{align}
    E^k \coloneqq \min\left\{ -g(x^{k}),\tfrac{1}{\nu_{k-1}}\mu^{k-1}\right\}, \label{eq:E}
\end{align}
which serves as a surrogate of approximate complimentary slackness (see  \cref{thm:E-KKT}).
In particular, we use this quantity to determine when the penalty parameter $\nu_k$ should be increased. 
\begin{algorithm}[t]
    \caption{\PBALM~((Inexact) proximal feasible-start ALM)}
    \label{alg:PBALM}
\begin{algorithmic}[1]
    \setlength\itemsep{0.5ex}
    \Require
        \begin{tabular}[t]{@{}l@{}}
            Initial feasible point \( x^0\in\dom f_2\) such that \(h( x^0)=0\), \(g(x^0) \leq 0\)
        \\
            Multipliers \(\lambda^0\in \R^p\) and $\mu^0\in\R_+^m$  and 
            penalties \(\rho_0, \nu_0, \gamma_0 \in \rp\)
         \\  
         Parameters \(\beta\in (0,1)\), \(\xi_1,\xi_2 \geq 1\), and \(\delta, \hat{\rho},  \hat{\nu},  \hat{\gamma} \in \rp\)
        \\
            Sequence of tolerances \((\tau_k)_{k\in\N} \subset \rp\) with \(\tau_k \to \tau\) for some \(\tau \in \R_+\)
        \\
            Functions \(\phi\), \(\phi_\gamma\) satisfying \cref{assump:penalty}
        \end{tabular}

    \item[For \(k=0,1,2\ldots\)]

    \State\label{state:proxALM:barx}
        \begin{tabular}[t]{@{}r@{~~}l@{}}
            Set& \(\hat{x}^k= x^k\) if \(\calL_{\rho_k, \nu_k, \gamma_k}(x^k,\lambda^k, \mu^k; x^k) + f_2(x^k) \le f( x^0) + \frac{1}{2\gamma_k}\|x^0 - x^k\|^2\),
            \\
            or& \(\hat{x}^k= x^0\) otherwise
        \end{tabular}
    \State\label{state:proxALM:x}
        \parbox[t]{0.98\linewidth}{
            Starting at \(\hat x^k\), find a point \( x^{k+1} \in \dom f_2\) approximately solving
            \begin{equation*}
                \minimize{}\left\{\calL_{\rho_k, \nu_k, \gamma_k}({}\cdot{},\lambda^k, \mu^k; x^k) + f_2({}\cdot{})\right\},
            \end{equation*}
            such that
            \begin{equation}\label{prop:proxALM:grad}
                \|x^{k+1}-\prox_{f_2}(x^{k+1}-\nabla_{x}\calL_{\rho_k, \nu_k, \gamma_k}( x^{k+1},\lambda^k, \mu^k; x^k))\|_\infty \leq \tau_k
            \end{equation}
        }
    
    \State\label{state:proxALM:lambda}
        Set \(\lambda^{k+1}=\lambda^k+\rho_kh( x^{k+1})\)
    
    \State\label{state:proxALM:mu}
        Set \(\mu^{k+1}=\max\{0,\mu^k+\nu_k g( x^{k+1})\}\)
    
    \State\label{state:proxALM:rho}
        \begin{tabular}[t]{@{}r@{~~}l@{}}
            Set & \(\rho_{k+1}=\rho_k\) ~if \(\|h( x^{k+1})\|_\infty\leq\beta\|h( x^k)\|_\infty\),\\
            or & \(\rho_{k+1}=\max\set{\xi_1\rho_k,\hat{\rho}\phi(k+1)}\) ~otherwise
        \end{tabular}
    
    \State\label{state:proxALM:nu}
        \begin{tabular}[t]{@{}r@{~~}l@{}}
            Set & \(\nu_{k+1}=\nu_k\) ~if \(\|E^{k+1}\|_\infty\leq\beta\|E^k\|_\infty\) where \(E_k\) is as in \eqref{eq:E},
            \\
            or & \(\nu_{k+1}=\max\set{\xi_2\nu_k,\hat{\nu}\phi(k+1)}\) ~otherwise
        \end{tabular}
        \State\label{state:proxALM:eta}
        Set $\gamk*=\max\{\delta\|x^0 - x^{k+1}\|^2, \hat{\gamma} \phi_\gamma(k+1)\}$
    \end{algorithmic}
\end{algorithm}
The main difference from the standard Augmented Lagrangian literature lies in the way the penalty parameters are updated. \PBALM{} introduces a general update rule using a function \(\phi\) that, while encompassing the standard linear growth, also allows for alternatives such as \(\phi(k) = (k+1)^\alpha\) for \(\alpha>1\), which grow more aggressively in the early stages but increases more moderately in the long runs. This approach helps mitigating the numerical issues  known to arise when a purely linear rate of increase is used. Our proposed update rule for the proximal parameter $\gamma_k$ uses a separate function $\phi_\gamma$:
\begin{align}
    \gamma_k = \max\left\{\delta\|x^0 - x^k\|^2, \gamma_0\phi_\gamma(k)\right\},\notag
\end{align}
where $\delta \in \rp$ is a constant. Our standing assumption for \(\phi\) and, separately, for $\phi_\gamma$ is as follows.
\begin{assumption}[Growth function]\label{assump:penalty}
    A function $\psi\colon\N\to\rp$ satisfies the growth condition if
\begin{align}
    \limsup_{k \to \infty} \tfrac{\psi(k+1)}{\psi(k)} < M, \quad   \lim_{k \to \infty}\tfrac{k}{\psi(k)}\to 0,\notag
\end{align}
for some $M \in \rp$. We require that both the function $\phi$, appearing in the update rule for $\rho_k$ and $\nu_k$, and the function $\phi_\gamma$, appearing in the update rule for $\gamma_k$, satisfy this growth condition.
\end{assumption}
In practice, we suggest to use a monotonically increasing \(\phi(k) = (k+1)^\alpha\), for some \(\alpha>1\), as it allows for a quick recovery from a small initialization of the penalty parameters, while it ameliorates ill-conditioning typically observed with standard linear increments. For the proximal parameter, we recommend choosing $\phi_\gamma(k) = (k+1)^{\alpha_\gamma}$ with $\alpha_\gamma>1$ close to $1$. We remark that $\rho_k$ and $\nu_k$ could in principle also use different growth functions from one another (each still satisfying \cref{assump:penalty}). For simplicity of presentation we keep them tied to the same function $\phi$, distinct from the function $\phi_\gamma$ used for the proximal parameter $\gamma_k$. We refer the reader to \cref{sec:setup} for detailed discussion on the choice of parameters.

At each iteration of \PBALM, a minimization problem involving the proximal augmented Lagrangian is solved. To ensure convergence, it is required that the algorithm provides a stationary point of the sum of the proximal augmented  Lagrangian and \(f_2\) with a maximum tolerance of $\tau_k \in \rp$ at iteration $k$, where $\seq{\tau_k}$ is a sequence of positive tolerances converging to a constant $\tau \in \R_+$.

Because \PBALM{} only requires $x^{k+1}$ to satisfy \eqref{prop:proxALM:grad} approximately, rather than to exactly solve \eqref{eq:proxALM-updates}, the inclusion \eqref{eq:proxALM-fixed-point} no longer holds verbatim for its iterates. Its second and third blocks remain exact regardless of how $x^{k+1}$ is computed, since $\tfrac{1}{\rho_k}(\lambda^k-\lambda^{k+1}) = -h(x^{k+1})$ and $\tfrac{1}{\nu_k}(\mu^k-\mu^{k+1}) \in -g(x^{k+1}) + \NN_{\R_+^m}(\mu^{k+1})$ follow directly from \cref{state:proxALM:lambda,state:proxALM:mu}. The first (primal stationarity) block, in contrast, picks up an additional error of size $\tau_k$, as made precise in \cref{thmi:proxALM:conv-step-1}.

The main convergence properties of \PBALM{} are stated in \cref{thm:ALM:convergence}.
\begin{theorem}[$\varepsilon$-KKT termination for \PBALM]\label{thm:ALM:convergence}
     Suppose that \cref{assump:prob,assump:penalty} hold. Consider the sequence $\seq{x^k,\lambda^k,\mu^k}$ generated by \PBALM{}. Then, the following hold:
     \begin{enumerate}
        \item \label{thm:prim-res}
        The sequence $\seq{\tfrac{1}{\gamk}\|x^{k+1}- x^k\|}[k\in \N]$ converges to zero;
        \item \label{thm:ALM:convergence:residual:vanishing} 
        There exists an infinite subsequence \(\mathcal K \subseteq \N\) along which $\seq{\tfrac{1}{\rhok^2}\|\lambda^{k+1} - \lambda^{k}\|^2+ \tfrac1{\nuk^2}\|\mu^{k+1} - \mu^{k}\|^2}[k \in \mathcal K]$ converges to zero;
         \item \label{thm:eps-kkt}
         An $\varepsilon$-KKT point is obtained in a finite number of iterations, for any $\varepsilon > \tau \in \rnn$.
     \end{enumerate}
\end{theorem}
The proof of \cref{thm:ALM:convergence} is deferred to \Secref{sec:convergence}, where several properties of \PBALM{} are derived. \cref{thm:prim-res} follows directly from \cref{thm:eta-vanish-1}. Observe that, as a direct implication of how $\lambda^k$ and $\mu^k$ are updated in \PBALM,
\begin{align}\label{eq:dual-residual}
    \tfrac{1}{\rhok^2}\|\lambda^{k+1} - \lambda^{k}\|^2 + \tfrac1{\nuk^2}\|\mu^{k+1} - \mu^{k}\|^2 = \|h(x^{k+1})\|^2 + \|E^{k+1}\|^2.
\end{align}
Hence, in \Secref{sec:convergence}, we prove \cref{thm:ALM:convergence:residual:vanishing} by showing that $\|h(x^{k+1})\|^2$ and $\|E^{k+1}\|^2$ tend to zero along an infinite subsequence $\mathcal K \subseteq \N$, i.e., \cref{thm:ALM:convergence:residual:vanishing} follows from \eqref{eq:dual-residual} together with \cref{thm:epsKKT_conv}. In view of \eqref{eq:proxALM-fixed-point}, and using the characterization of an $\varepsilon$-KKT point that follows from \cref{thm:E-KKT}, we consequentially prove \cref{thm:eps-kkt}.

As shown in \Cref{thm:prim-res}, the quantity $\tfrac1{\gamma_k}\|x^{k+1}-x^k\|$ vanishes independently of the tolerances $\tau_k$. This follows from $\gamma_k$ growing faster than linearly (\cref{assump:penalty}) while $\tfrac1{2\rho_k}\|\lambda^k\|^2+\tfrac1{2\nu_k}\|\mu^k\|^2$ grows at most linearly (\cref{thm:proxALM:multipliers_bound}), as shown in the proof of \cref{thm:eta-vanish}. As noted above, this is precisely the correction that the first block of \eqref{eq:proxALM-fixed-point} acquires once $x^{k+1}$ is only computed inexactly. By \cref{thmi:proxALM:conv-step-1}, the genuine stationarity gap of $(x^{k+1},\lambda^{k+1},\mu^{k+1})$ is bounded by $\tau_k$ plus this vanishing quantity, and therefore only gets arbitrarily close to $\tau \coloneqq \lim_k \tau_k$. Hence, \cref{thm:eps-kkt} guarantees an $\varepsilon$-KKT point for any $\varepsilon>\tau$, rather than as $\varepsilon\to0$.

\begin{remark}[AKKT/AM-stationarity]\label{rem:akkt}
    Combined with \cref{thm:ALM:convergence:residual:vanishing}, the above shows that a selection of the KKT operator $T$ at $(x^{k+1},\lambda^{k+1},\mu^{k+1})$, as in \eqref{eq:KKT}, vanishes along $\mathcal K$ up to the tolerance $\tau_k$. This is the sequential residual condition defining an approximate KKT (AKKT) point \cite{andreani2010new,andreani2016cone}, extended to composite problems of the form \eqref{eq:prob} as AM-stationarity in \cite{de2023constrained}, except that the residual is driven to $\tau$ rather than to zero, so that the two coincide exactly when $\tau=0$. {Accordingly, i}f $\tau=0$ and $x^k \to x^\star$ along a subsequence of $\mathcal K$, then $x^\star$ is AM-stationary, and a sequential constraint qualification at $x^\star$ upgrades this to a genuine KKT point. Our focus is instead the finite-time attainment of an $\varepsilon$-KKT point for a fixed tolerance $\varepsilon>\tau$.
\end{remark}

\begin{remark}[Penalty parameters and initial feasibility]\label{rem:matrix-stepsize}
    \begin{enumerate}
    \item
    The use of scalar penalty parameters \(\rhok, \nuk\) for equality and inequality constraints in \PBALM{} is only for clarity of exposition, and our theoretical results extend trivially to the case of (block) diagonal penalty matrices. Specifically, given any decomposition \(h(x) = (h_1(x), h_2(x), \ldots, h_{n_h}(x))\), we may asign \(\rho_{i,k}\) as the penalty parameter corresponding to \(h_i\) with constraint violation measured in terms of \(\|h_i(x^k)\|_\infty\).

    \item 
        In many applications, such as those arising in optimal control, an initial feasible point can be obtained at little to no cost by simply unrolling the dynamics through a forward pass. More generally, however, a feasible point may be computed by solving a so-called \emph{phase I} problem, as is common in the literature on barrier-type methods \cite{conn2000trust}, see also \cite[Section 11.4]{boyd2004convex}. A brief discussion of this approach, in which \PBALM{} itself may be used to solve the phase I problem is given in \appref{sec:phase-I}.
        
        \item 
        \PBALM{} uses the same initial point \(x^0\) in \cref{state:proxALM:barx} throughout iterations. However, in practice, should the algorithm stumble upon a new feasible point along the iterates, the new feasible point may be used to replace \(x^0\).
    \end{enumerate}
\end{remark}

\section{Convergence analysis}\label{sec:convergence}
We begin by presenting a series of results on the properties of the sequence generated by \PBALM{}. These results are instrumental for the proof of the main result in \cref{thm:ALM:convergence} at the end of this section. 

First, we present a simple observation in \cref{thm:proxALM:L-bound} that the augmented Lagrangian is upper bounded by a quantity that depends solely on the variable at the previous iterate and the initial feasible point. The lemma is a trivial consequence of the initialization for the Lagrangian subproblem \eqref{prop:proxALM:grad} in \cref{state:proxALM:barx} of \PBALM{}, provided that \cref{state:proxALM:x} is solved by a descent method.
\begin{lemma}[Lagrangian upper bound]\label{thm:proxALM:L-bound}
    Let \cref{assump:feas_sol} hold, and let $x^0$ be an initial primal point of \PBALM{} such that $h(x^0) = 0$ and $g(x^0) \leq 0$. Suppose further that $x^{k+1}$ in \cref{state:proxALM:x} of \PBALM{} is computed by a descent method, so that
    \begin{align*}
        \calL_{\rho_k,\nu_k,\gamma_k}(x^{k+1},\lambda^k,\mu^k;x^k) + f_2(x^{k+1}) \leq \calL_{\rho_k,\nu_k,\gamma_k}(\hat{x}^k,\lambda^k,\mu^k;x^k) + f_2(\hat{x}^k).
    \end{align*}
    Then, for any $k\in \N$,
    \begin{align}
        \calL_{\rho_k,\nu_k,\gamma_k}(x^{k+1},\lambda^k,\mu^k;x^k) + f_2(x^{k+1}) \leq f(x^0) + \tfrac{1}{2\gamma_k}\|x^0 - x^k\|^2.\label{eq:proxALM:L-bound}
    \end{align}
\end{lemma}
We next show that the scaled sum of the multipliers, \(\tfrac{1}{2\rhok*}\|\lambda^{k+1}\|^2 + \tfrac{1}{2\nuk*}\|\mu^{k+1}\|^2\), can increase by at most a fixed constant from one iteration to the next. 
This key result implies, in particular, that the scaled primal residual \(\tfrac{1}{\gamma_k}\|x^{k+1} - x^k\|\) vanishes asymptotically.

\begin{theorem}[Multipliers control and primal residual]\label{thm:eta-vanish}
    Suppose that \cref{assump:prob,assump:penalty} hold, and consider the iterates generated by \PBALM{}. Then, the following hold:
    \begin{enumerate}
        \item \label{thm:proxALM:multipliers_bound} For any \(k\in \N\)
        \begin{align*}
            \tfrac{1}{2\rhok*}\|\lambda^{k+1}\|^2 + \tfrac{1}{2\nuk*}\|\mu^{k+1}\|^2 
            \leq \tfrac{1}{2\rho_k}\|\lambda^k\|^2 + \tfrac{1}{2\nu_k}\|\mu^k\|^2
            +  c_1 - \tfrac{1}{2\gamma_k}\|x^{k+1} - x^k\|^2,
        \end{align*}
        and in particular,
        \begin{align}
            \tfrac{1}{2\rhok}\|\lambda^{k}\|^2 + \tfrac{1}{2\nuk}\|\mu^{k}\|^2 \leq 
            \tfrac{1}{2\rho_0}\|\lambda^0\|^2 + \tfrac{1}{2\nu_0}\|\mu^0\|^2
            +  kc_1,
            \label{eq:proxALM:multipliers_bound1}
        \end{align}
        where $c_1 \coloneqq f(x^0) - f^\star + \tfrac{1}{2\delta}$;
        \item \label{thm:eta-vanish-1} $\tfrac{1}{\gamma_k}\|x^{k+1} - x^k\| \to 0 \text{ as } k \to \infty$.
    \end{enumerate} 
\end{theorem}
\begin{proof}
    We can rewrite~\eqref{eq:prox-ALM:lagrangian} as
    \begin{align}
        \calL_{\rho,\nu,\gamma}(x,\lambda, \mu;v) &= f_1(x) + \tfrac{1}{2\rho}\|\rho h(x) + \lambda\|^2 + \tfrac{1}{2\nu}\|\max\{0, \nu g(x) + \mu\}\|^2 - \tfrac{1}{2\rho}\|\lambda\|^2 - \tfrac{1}{2\nu}\|\mu\|^2\notag\\
        &\qquad + \tfrac{1}{2\gamma}\|x - v\|^2.\label{eq:ALM:proof:2}
    \end{align}
    Using~\cref{state:proxALM:lambda,state:proxALM:mu} in~\eqref{eq:ALM:proof:2} gives
    \begin{align*}
        \tfrac{1}{2\rho_k}\|\lambda^{k+1}\|^2 + \tfrac{1}{2\nu_k}\|\mu^{k+1}\|^2 &= \calL_{\rho_{k},\nu_{k},\gamma_k}(x^{k+1}, \lambda^k, \mu^k;x^k) - f_1(x^{k+1}) - \tfrac{1}{2\gamma_k}\|x^{k+1} - x^k\|^2\\
        &\quad + \tfrac{1}{2\rho_k}\|\lambda^k\|^2 + \tfrac{1}{2\nu_k}\|\mu^k\|^2\\
        &\le f(x^0) - f^\star + \tfrac{1}{2\gamma_{k}}\|x^0 - x^{k}\|^2 - \tfrac{1}{2\gamma_k}\|x^{k+1} - x^k\|^2\\
        &\quad + \tfrac{1}{2\rho_k}\|\lambda^k\|^2 + \tfrac{1}{2\nu_k}\|\mu^k\|^2,
    \end{align*}
    where the inequality follows from \cref{thm:proxALM:L-bound} and \cref{assump:obj}. Now, using the fact that $\gamma_k \ge \delta \|x^0 - x^k\|^2$ as ensured by \cref{state:proxALM:eta} of \PBALM, we have
    \begin{align*}
        \tfrac{1}{2\rho_k}\|\lambda^{k+1}\|^2 + \tfrac{1}{2\nu_k}\|\mu^{k+1}\|^2 &\le c_1 - \tfrac{1}{2\gamma_k}\|x^{k+1} - x^k\|^2 + \tfrac{1}{2\rho_k}\|\lambda^k\|^2 + \tfrac{1}{2\nu_k}\|\mu^k\|^2.
    \end{align*}
     Since $\rho_k$ and $\nu_k$ are monotonically nondecreasing, this implies that
     \begin{align}
        \tfrac{1}{2\rho_{k+1}}\|\lambda^{k+1}\|^2 + \tfrac{1}{2\nu_{k+1}}\|\mu^{k+1}\|^2 &\le \tfrac{1}{2\rho_k}\|\lambda^{k+1}\|^2 + \tfrac{1}{2\nu_k}\|\mu^{k+1}\|^2\notag\\
        &\le c_1 - \tfrac{1}{2\gamma_k}\|x^{k+1} - x^k\|^2 + \tfrac{1}{2\rho_k}\|\lambda^k\|^2 + \tfrac{1}{2\nu_k}\|\mu^k\|^2\label{eq:proxALM:multipliers_bound0}
        \\ 
        &\le \tfrac{1}{2\rho_k}\|\lambda^k\|^2 + \tfrac{1}{2\nu_k}\|\mu^k\|^2
        + c_1. 
        \nonumber
    \end{align}
    A telescoping argument then yields the bound stated in \eqref{eq:proxALM:multipliers_bound1}. 
    Now, rearranging \eqref{eq:proxALM:multipliers_bound0} gives
    \begin{align*}
        \tfrac{1}{2\gamma_k^2}\|x^{k+1} - x^k\|^2 &\le \tfrac{1}{\gamma_k}\left(c_1 + \tfrac{1}{2\rho_k}\|\lambda^k\|^2 + \tfrac{1}{2\nu_k}\|\mu^k\|^2\right)\\
        &\le \tfrac{1}{\gamma_k}\left((k+1)c_1 + \tfrac{1}{2\rho_0}\|\lambda^0\|^2 + \tfrac{1}{2\nu_0}\|\mu^0\|^2\right)\\
        &\le \tfrac{1}{\hat{\gamma}\phi_\gamma(k)}\left((k+1)c_1 + \tfrac{1}{2\rho_0}\|\lambda^0\|^2 + \tfrac{1}{2\nu_0}\|\mu^0\|^2\right),
    \end{align*}
    where the second inequality is due to \eqref{eq:proxALM:multipliers_bound1} and the third inequlity uses the fact that $\gamma_k \ge \hat{\gamma}\phi_\gamma(k)$ (see \cref{state:proxALM:eta} of \PBALM). Consequently, using \cref{assump:penalty} applied to $\phi_\gamma$, we find $\tfrac{1}{\gamma_k}\|x^{k+1} - x^k\| \to 0$ as $k \to \infty$.
\end{proof}
\begin{lemma}[Fixed-point residual]\label{thm:proxALM:conv-step-1}
    Suppose that \cref{assump:prob,assump:penalty} hold, and let \(\seq{x^k, \lambda^k, \mu^k}\) be the sequence of iterates generated by \PBALM. Then, the following properties hold:
    \begin{enumerate}
        \item
        \(
            \nabla_x \calL(x^{k+1},\lambda^{k+1},\mu^{k+1}) = \nabla_x \calL_{\rho_{k},\nu_{k},\gamma_{k}}(x^{k+1},\lambda^k,\mu^k;x^k) - \tfrac{1}{\gamma_k}(x^{k+1} - x^k);
        \)
        \item \label{thmi:proxALM:conv-step-1}
        \(\|x^{k+1} - \prox_{f_2}(x^{k+1} - \nabla_{x}\calL(x^{k+1},\lambda^{k+1},\mu^{k+1}))\|_\infty \leq \varepsilon,
        \) for $k$ large enough and any $\varepsilon > \tau$.
    \end{enumerate} 
\end{lemma}
\begin{proof}
    Using~\cref{state:proxALM:lambda,state:proxALM:mu} in \eqref{eq:L}, we get that
    \begin{align}
        \calL(x,\lambda^{k+1},\mu^{k+1}) &= \calL_{\rho_{k},\nu_{k},\gamma_{k}}(x,\lambda^{k},\mu^{k}; x^k) - \tfrac{\rho_{k}}{2}\|h(x) - h(x^{k+1})\|^2 + \tfrac{\rho_k}{2}\|h(x^{k+1})\|^2 + \tfrac{1}{2\nu_{k}}\|\mu^{k}\|^2\nonumber\\
        & \quad + \langle\max\{0,\mu^{k}+\nu_{k}g(x^{k+1})\},g(x)\rangle - \tfrac{1}{2\nu_{k}}\|\max\{0,\mu^{k}+\nu_{k}g(x)\}\|^2\nonumber\\
        & \quad - \tfrac{1}{2\gamma_{k}}\|x - x^k\|^2,\notag
    \end{align}
    and
    \begin{align*}
        \nabla_x \calL(x,\lambda^{k+1},\mu^{k+1}) &= \nabla_x \calL_{\rho_{k},\nu_{k},\gamma_{k}}(x,\lambda^k,\mu^k;x^k) - \rho_k J_h(x)^\top(h(x) - h(x^{k+1}))\\
        & \quad + J_g(x)^\top(\max\{0,\mu^{k}+\nu_{k}g(x^{k+1})\}) - J_g(x)^\top(\max\{0,\mu^{k}+\nu_{k}g(x)\})\\
        & \quad - \tfrac{1}{\gamma_k}(x - x^k).
    \end{align*}
    Hence, we have
    \begin{align}
        \nabla_x \calL(x^{k+1},\lambda^{k+1},\mu^{k+1}) = \nabla_x \calL_{\rho_{k},\nu_{k},\gamma_{k}}(x^{k+1},\lambda^k,\mu^k;x^k) - \tfrac{1}{\gamma_k}(x^{k+1} - x^k).\label{eq:proxALM:conv-lemma:eq1}
    \end{align}
    By the triangle inequality, we have
\begin{align*}
    &\|x^{k+1} - \prox_{f_2}(x^{k+1} - \nabla_x \calL(x^{k+1},\lambda^{k+1},\mu^{k+1}))\|_\infty \\
    &\le~ \|x^{k+1} - \prox_{f_2}(x^{k+1} - \nabla_x \calL(x^{k+1},\lambda^{k+1},\mu^{k+1}) - \tfrac{1}{\gamma_k}(x^{k+1} - x^k))\|_\infty \\
    &\quad + \|\prox_{f_2}(x^{k+1} - \nabla_x \calL(x^{k+1},\lambda^{k+1},\mu^{k+1}) - \tfrac{1}{\gamma_k}(x^{k+1} - x^k)) \\
    &\qquad - \prox_{f_2}(x^{k+1} - \nabla_x \calL(x^{k+1},\lambda^{k+1},\mu^{k+1}))\|_\infty \\
    &\le~ \|x^{k+1} - \prox_{f_2}(x^{k+1} - \nabla_x \calL_{\rho_{k},\nu_{k},\gamma_{k}}(x^{k+1},\lambda^k,\mu^k;x^k))\|_\infty 
    + \tfrac{1}{\gamma_k}\|x^{k+1} - x^k\|,
\end{align*}
where the second inequality results from \eqref{eq:proxALM:conv-lemma:eq1}, the nonexpansiveness of the proximal operator and $\|\cdot\|_\infty \le \|\cdot\|$. Now, by \eqref{prop:proxALM:grad} of \cref{state:proxALM:x} in \PBALM, we get
    \begin{align*}
        &\|x^{k+1} - \prox_{f_2}(x^{k+1} - \nabla_x \calL(x^{k+1},\lambda^{k+1},\mu^{k+1}))\|_\infty \le \tau_k + \tfrac{1}{\gamma_k}\|x^{k+1} - x^k\|.
    \end{align*}
    Since $\tau_k \to \tau$ and, by \cref{thm:eta-vanish-1}, $\tfrac{1}{\gamma_k}\|x^{k+1}-x^k\| \to 0$, for $k$ large enough the right-hand side is smaller than any $\varepsilon > \tau$.
\end{proof}
For every iteration $k$ at which the penalty parameters increase in \cref{state:ALM:rho,state:ALM:nu} of \PBALM, we derive upper bounds on the feasibility and complementarity of $(x^k, \lambda^k, \mu^k)$ in the following lemma.
\begin{lemma}[Feasibility and complementarity control]\label{lem:feasibility:bound}
    Suppose that \cref{assump:prob,assump:penalty} hold, and consider the sequence \(\seq{x^k, \lambda^k, \mu^k}\) of iterates generated by \PBALM. Define the index sets
    \begin{align}
        \calK_\rho \coloneqq \{k \in \N  \mid \|h(x^k)\|_\infty > \beta \|h(x^{k-1})\|_\infty\}, \quad \calK_\nu \coloneqq \{k \in \N  \mid \|E^k\|_\infty > \beta \|E^{k-1}\|_\infty\},
        \label{eq:indexsets}
    \end{align}
    and let
        \begin{align*}
            A_k \coloneqq 2k(f(x^0) - f^\star + \tfrac{1}{2\delta}) + \tfrac{1}{\rho_0}\|\lambda^0\|^2 + \tfrac{1}{\nu_0}\|\mu^0\|^2.
        \end{align*}
    Then, the following hold:
    \begin{enumerate}
        \item \label{lem:Ek:bound}For any \(k \in \calK_\nu\),
        \begin{align*}
            \tfrac{1}{2}\|E^{k+1}\|^2 &\le \tfrac{1}{\hat{\nu}}\max\left\{\tfrac{\xi_2}{\phi(k)}, \tfrac{\phi(k+1)}{\phi(k)^2}\right\}A_{k+1} + \tfrac{1}{\hat{\nu}\phi(k)}A_k; 
            \numberthis\label{eq:case2:Ekto0}
        \end{align*} 
        \item \label{lem:hk:bound} For any \(k \in \calK_\rho\),
            \begin{align*}
            \tfrac{1}{2}\|h(x^{k+1})\|^2 &\leq \tfrac{1}{ \hat{\rho}}\max\left\{\tfrac{\xi_1}{\phi(k)}, \tfrac{\phi(k+1)}{\phi(k)^2}\right\}A_{k+1} + \tfrac{1}{\hat{\rho}\phi(k)}A_k. 
            \numberthis\label{eq:case2:hkto0}
        \end{align*} 
    \end{enumerate}
\end{lemma}
\begin{proof}
Observe that
\begin{align}
    \|E^{k+1}\|^2 &= \|\min\{-g(x^{k+1}),\tfrac{1}{\nu_k}\mu^k\}\|^2\notag\\
    &\le 2\|\max\{g(x^{k+1}),-\tfrac{1}{\nu_k}\mu^k\} + \tfrac{1}{\nu_k}\mu^k\|^2 + 2\|\tfrac{1}{\nu_k}\mu^k\|^2\nonumber\\
    &= \tfrac2{\nuk^2}\|\mu^{k+1}\|^2 + \tfrac2{\nuk^2}\|\mu^k\|^2.\label{eq:ALM:cons-bounds-1}
\end{align}
It follows from \eqref{eq:ALM:cons-bounds-1} and \cref{thm:proxALM:multipliers_bound} that 
\begin{align*}
    \tfrac{1}{2}\|E^{k+1}\|^2 &\le \tfrac{1}{\nu_k^2}(\|\mu^{k+1}\|^2 + \|\mu^k\|^2)\\
    &\le \tfrac{\nu_{k+1}}{\nu_k^2}A_{k+1} + \tfrac{1}{\nu_k}A_k\\
    &\leq \max\left\{\tfrac{\xi_2}{\nu_k}, \tfrac{\hat{\nu}\phi(k+1)}{\nu_k^2}\right\}A_{k+1} + \tfrac{1}{\nu_k}A_k,
\end{align*}
where the last inequality follows from the fact that, regardless of the condition in \cref{state:proxALM:nu} of \PBALM, \(\nuk* \leq \max\{\xi_2 \nuk, \hat{\nu}\phi(k+1)\}\) since \(\xi_2 \geq 1\).
For any \(k \in \calK_\nu\) it holds that \(\nu_k \geq \hat{\nu} \phi(k)\). The claimed inequality \eqref{eq:case2:Ekto0} follows immediately by using this lower bound for \(\nuk\) in the above inequality. 

The second claim follows similarly, this time noting that 
\begin{align}
\|h(x^{k+1})\|^2 \le{}& 2\|h(x^{k+1}) + \tfrac{1}{\rho_k}\lambda^k\|^2 + 2\|\tfrac{1}{\rho_k}\lambda^k\|^2
\nonumber
\\ 
={}& \tfrac2{\rhok^2}\|\lambda^{k+1}\|^2 + \tfrac2{\rhok^2}\|\lambda^k\|^2
\notag
\end{align}
and using the bound \(\rhok \geq \hat{\rho}\phi(k)\). 
\end{proof}

We next obtain the following result which is used to establish that the sequence of iterates of \PBALM{} is such that the feasibility and complementarity residuals tend to zero along an infinite subsequence.
\begin{theorem}[Feasibility and complementarity residuals]\label{thm:epsKKT_conv}
    Suppose that \cref{assump:prob,assump:penalty} hold, and let \(\seq{x^k, \lambda^k, \mu^k}\) be the sequence of iterates generated by \PBALM. Then, for any $\varepsilon>0$, there exists an infinite subsequence $\mathcal K \subseteq \N$ such that, for all \(k \in \mathcal K\), it holds that \(\|h(x^k)\|_\infty \leq \varepsilon\) and \(\|E^k\|_\infty \leq \varepsilon\).
\end{theorem}
\begin{proof}

For $k\in\N$ define
\[
    t_\rho(k) \coloneqq \max\set{i \le k}[i \in \calK_\rho],
    \qquad
    t_\nu(k) \coloneqq \max\set{i \le k}[i \in \calK_\nu],
\]
with the convention $\max\emptyset \coloneqq 0$, and let
\[
    \calS_\rho \coloneqq \set{j \in \calK_\rho}[j-1 \notin \calK_\rho],
    \qquad
    \calS_\nu \coloneqq \set{j \in \calK_\nu}[j-1 \notin \calK_\nu]
\]
denote the indices at which a run of consecutive increases begins. 
Since $\rho_k$ is nondecreasing and equals $\max\{\xi_1\rho_{k-1}, \hat{\rho}\phi(k)\} \ge \hat\rho\phi(k)$ at every index at which it increases, and remains constant otherwise,
\begin{equation}
    \label{eq:pen-lower}
    \rho_k \ge \hat{\rho}\,\phi(t_\rho(k))
    \quad\text{whenever }\calK_\rho \cap [0,k] \neq \emptyset,
\end{equation}
and analogously $\nu_k \ge \hat{\nu}\phi(t_\nu(k))$ whenever $\calK_\nu \cap [0,k] \neq \emptyset$.
We proceed with intermediate claims.

\emph{Claim 1: $\|h(x^j)\| \to 0$ as $\N \setminus \calS_\rho \ni j \to \infty$, and $\|E^j\| \to 0$ as $\N \setminus \calS_\nu \ni j \to \infty$.}

\noindent{}Let $j \notin \calS_\rho$. If $j-1 \in \calK_\rho$, then \cref{lem:hk:bound} applied at $j-1$ together with \cref{assump:penalty} gives $\|h(x^{j})\| \to 0$. Otherwise $j - 1 \notin \calK_\rho$, and since $j \notin \calS_\rho$ this forces $j \notin \calK_\rho$; setting $t \coloneqq t_\rho(j) < j$, the indices $t+1,\dots,j$ all lie outside $\calK_\rho$, so that $j-t-1$ applications of \eqref{eq:indexsets} yield
\[
    \|h(x^{j})\|_\infty \le \beta\|h(x^{j-1})\|_\infty \le \ldots \le \beta^{\,j-t-1}\|h(x^{t+1})\|_\infty .
\]
Two cases arise. If $\calK_\rho$ is infinite, then $t = t_\rho(j) \to \infty$ as $j \to \infty$, so that $\|h(x^{t+1})\| \to 0$ by \cref{lem:hk:bound} applied at $t \in \calK_\rho$, while the factor $\beta^{\,j-t-1}$ is bounded by $1$ since $\beta \in (0,1)$ and $j \ge t+1$. If instead $\calK_\rho$ is finite, then $t_\rho(j)$ is constant and equal to $\bar t \coloneqq \max \calK_\rho$ for every $j > \bar t$, so that $\|h(x^{\bar t+1})\|_\infty$ is a fixed constant while the factor $\beta^{\,j-\bar t-1}$ tends to zero as $j \to \infty$. In either case the right-hand side tends to zero. (If $\calK_\rho = \emptyset$, the same argument applies with $\bar t = 0$.) Thus,  $\|h(x^{j})\| \to 0$. The argument for $E^j$ is identical.

\emph{Claim 2: if $\N \setminus (\calS_\rho \cup \calS_\nu)$ is infinite, the assertion of the theorem holds.} 

\noindent{}This set is contained in both $\N \setminus \calS_\rho$ and $\N \setminus \calS_\nu$, so by Claim 1 both $\|h(x^j)\|$ and $\|E^j\|$ tend to zero along it; any $\varepsilon > 0$ is thus attained on an infinite subset of it.

\emph{Claim 3: if $\N \setminus (\calS_\rho \cup \calS_\nu)$ is finite, the assertion of the theorem holds.} 

\noindent{}
Suppose first that $\calS_\rho \cap \calS_\nu$ is infinite. Then \cref{lem:hk:bound,lem:Ek:bound} both apply at every $j \in \calS_\rho \cap \calS_\nu \subseteq \calK_\rho \cap \calK_\nu$, and by \cref{assump:penalty}, the assertion holds along the infinite set $\set{j+1}[j \in \calS_\rho \cap \calS_\nu]$.

Suppose instead that $\calS_\rho \cap \calS_\nu$ is finite. 
By definition neither $\calS_\rho$ nor $\calS_\nu$ contains two consecutive integers. Let $J \in \N$ be such that $\calS_\rho \cup \calS_\nu \supseteq \set{j \in \N}[j \ge J]$.
Enlarging $J$ if necessary, every $j \ge J$ then belongs to exactly one of $\calS_\rho$ and $\calS_\nu$, and since neither set contains two consecutive integers, membership alternates. Consequently, for every $j \ge J+1$, one of the two consecutive indices $j-1, j$ belongs to $\calK_\rho$ and the other to $\calK_\nu$, so that
\[
    t_\rho(j) \ge j-1
    \qquad\text{and}\qquad
    t_\nu(j) \ge j-1 .
\]
By \cref{assump:penalty} there exists $M > 1$ with $\phi(k+1) < M\phi(k)$ for all $k$ large. Hence $\phi(j-1) > \phi(j)/M$, and since $t_\rho(j) \in \{j-1, j\}$ we obtain $\phi(t_\rho(j)) \ge \phi(j)/M$ in either case. The same holds for $t_\nu(j)$. Combining with \eqref{eq:pen-lower},
\[
    \rho_j \ge \tfrac{\hat\rho}{M}\phi(j),
    \qquad
    \nu_j \ge \tfrac{\hat\nu}{M}\phi(j),
\]
repeating the derivation of \eqref{eq:case2:hkto0} with $\rho_j \ge \tfrac{\hat\rho}{M}\phi(j)$ in place of $\rho_j \ge \hat\rho\phi(j)$, together with $\phi(j+1) < M\phi(j)$, gives
\[
    \tfrac{1}{2}\|h(x^{j+1})\|^2
    \le \tfrac{M}{\hat{\rho}}\max\left\{\tfrac{\xi_1}{\phi(j)}, \tfrac{M^2}{\phi(j)}\right\}A_{j+1} + \tfrac{M}{\hat{\rho}\phi(j)}A_j
    = \tfrac{M\max\{\xi_1, M^2\}}{\hat{\rho}}\tfrac{A_{j+1}}{\phi(j)} + \tfrac{M}{\hat{\rho}}\tfrac{A_j}{\phi(j)} ,
\]
and analogously for $\tfrac{1}{2}\|E^{j+1}\|^2$, repeating the derivation of \eqref{eq:case2:Ekto0} with $\nu_j \ge \tfrac{\hat\nu}{M}\phi(j)$ in place of $\nu_j \ge \hat\nu\phi(j)$. Since $A_k = O(k)$ and $k/\phi(k) \to 0$ by \cref{assump:penalty}, both right-hand sides tend to zero, so both residuals tend to zero at the common index $j+1$ for every large $j$.

Claims 2 and 3 are exhaustive, which concludes the proof.
\end{proof}

Combining the above results, we can conclude that \PBALM{} finds an approximate KKT point after a finite number of iterations. To formalize the proof statement about this convergence property (\ie, as stated in \cref{thm:ALM:convergence}), we state the following instrumental result about the quantity $E^k$ defined in \eqref{eq:E} for quantifying inequality constraint violations and complementarity. The proof is given in \appref{app:supp_proofs} for completeness.
\begin{lemma}[\protect{\cite[Lemma 3.1]{birgin2020complexity}}]\label{thm:E-KKT}
    Suppose that \cref{assump:feas_sol} holds, and let $E^k$ be defined as in~\eqref{eq:E}. Then, for all $k\in \N$ and any $\varepsilon \in \rnn$, $\|E^{k+1}\|_\infty \leq \varepsilon$ implies that $\|\max\{0,g(x^{k+1})\}\|_\infty \leq \varepsilon$ and $\mu_i^{k+1} = 0$ for all $i \in [m]$ whenever $g_i(x^{k+1}) < -\varepsilon$.
\end{lemma}
We are now ready to combine all the above results to establish that the algorithm finds an approximate KKT point in finite number of iterations.
\begin{appendixproof}{thm:ALM:convergence}
    \cref{thm:prim-res} follows directly from \cref{thm:eta-vanish-1}, and \cref{thm:ALM:convergence:residual:vanishing} follows from \eqref{eq:dual-residual} together with \cref{thm:epsKKT_conv}.
    As shown in \cref{thmi:proxALM:conv-step-1}, it holds that for \(k\) large enough $\|x^k - \prox_{f_2}(x^k - \nabla_{x}\calL(x^k,\lambda^k,\mu^k))\|_\infty \leq \varepsilon$. Moreover, by \cref{thm:epsKKT_conv}, there exists a subsequence along which both $\|h(x^k)\|_\infty \leq \varepsilon$ and $\|E^k\|_\infty \leq \varepsilon$ hold. Now, note that by \cref{thm:E-KKT} and \cref{def:epsKKT}, these imply that $(x^k,\lambda^k,\mu^k)$ is an $\varepsilon$-KKT point of problem~\eqref{eq:prob}.
    Consequently, we conclude that an $\varepsilon$-KKT point of \eqref{eq:prob} is obtained after a finite number of iterations.
\end{appendixproof}
\section{An inexact ALM as a limiting case}\label{sec:ALM-algorithm}
\begin{algorithm}[t]
	\caption{\BALM~((Inexact) feasible-start ALM)}
	\label{alg:BALM}
\begin{algorithmic}[1]
    \setlength\itemsep{0.5ex}
    \Require
        \begin{tabular}[t]{@{}l@{}}
            Initial feasible point \( x^0\in\dom f_2\) such that \(h( x^0)=0\), \(g(x^0) \leq 0\)
        \\
            Multipliers \(\lambda^0\in \R^p\) and $\mu^0\in\R_+^m$  and 
            penalties \(\rho_0, \nu_0 \in \rp\)
        \\
            Parameters \(\beta \in (0,1)\), \(\xi_1,\xi_2 \geq 1\), \(\hat{\rho}, \hat{\nu} \in \rnn\)
        \\
            Sequence of tolerances \((\tau_k)_{k\in\N} \subset \rp\) with \(\tau_k \to \tau\) for some \(\tau \in \R_+\)
        \\
            Function \(\phi\) satisfying \cref{assump:penalty} {\small (e.g., \(\phi(k) = (k+1)^\alpha\) for some \(\alpha>1\))}
        \end{tabular}

    \item[For \(k=0,1,2\ldots\)]
    \State\label{state:ALM:barx}
        Set \(\hat{x}^k= x^k\) if \(\calL_{\rho_k, \nu_k}(x^k,\lambda^k, \mu^k) + f_2(x^k) \le f( x^0)\), or \(\hat{x}^k= x^0\) otherwise

    \State\label{state:ALM:x}
        \parbox[t]{0.98\linewidth}{
            Starting at \(\hat x^k\), find a point \( x^{k+1} \in \dom f_2\) approximately solving
            \begin{equation*}
                \minimize{}\left\{\calL_{\rho_k, \nu_k}({}\cdot{},\lambda^k, \mu^k) + f_2({}\cdot{})\right\},
            \end{equation*}
            such that
            \begin{equation*}
                \|x^{k+1}-\prox_{f_2}( x^{k+1}-\nabla_{x}\calL_{\rho_k, \nu_k}( x^{k+1},\lambda^k, \mu^k))\|_\infty \leq\tau_k
            \end{equation*}
        }

    \State\label{state:ALM:lambda}
        Set \(\lambda^{k+1}=\lambda^k+\rho_kh( x^{k+1})\)

    \State\label{state:ALM:mu}
        Set \(\mu^{k+1}=\max\{0,\mu^k+\nu_k g( x^{k+1})\}\)
    \State\label{state:ALM:rho}
        \begin{tabular}[t]{@{}r@{~~}l@{}}
            Set & \(\rho_{k+1}=\rho_k\) ~if \(\|h( x^{k+1})\|_\infty\leq\beta\|h( x^k)\|_\infty\),\\
            or & \(\rho_{k+1}=\max\set{\xi_1\rho_k,\hat{\rho}\phi(k+1)}\) ~otherwise 
        \end{tabular}
    
    \State\label{state:ALM:nu}
        \begin{tabular}[t]{@{}r@{~~}l@{}}
            Set & \(\nu_{k+1}=\nu_k\) ~if \(\|E^{k+1}\|_\infty\leq\beta\|E^k\|_\infty\), where 
            \(E^k = \min\left\{ -g(x^{k}),\tfrac{1}{\nu_{k-1}}\mu^{k-1}\right\} \)
            \\
            or & \(\nu_{k+1}=\max\set{\xi_2\nu_k,\hat{\nu}\phi(k+1)}\) ~otherwise
        \end{tabular}
    \end{algorithmic}
    
\end{algorithm}
The results obtained for \PBALM{} in the previous section extend to its ALM counterpart which does not involve a proximal term. In this section, we present \refBALM, an ALM version of \PBALM, and briefly discuss how some results on \PBALM, particularly \cref{thm:ALM:convergence}, hold also for \BALM~as special cases. For ease of our presentation, we define
\begin{align}
    \calL_{\rho, \nu}(x,\lambda,\mu) \coloneqq f_1(x) + \langle \lambda, h(x) \rangle + \frac{\rho}{2} \|h(x)\|^2 + \frac{1}{2\nu} \|\max\left\{0, \nu g(x) + \mu\right\}\|^2 - \frac1{2\nu}\|\mu\|^2, \label{eq:L'}
\end{align}
which is the augmented Lagrangian for problem \eqref{eq:prob} without the proximal term.

As a first instance, the following is analogous, for \BALM, to the feasibility property \eqref{eq:proxALM:L-bound} (via a similar argument):
\begin{align}
    \calL_{\rho_k,\nu_k}(x^{k+1},\lambda^k,\mu^k) + f_2(x^{k+1}) \leq f(x^0),\label{eq:ALM:L-bound}
\end{align}
for every $k\in \N$. This property of \BALM~helps to initialize subproblems properly, as well as establish the following result, whose analogue for \PBALM{} appears in the proof of \cref{thm:eta-vanish}.
\begin{lemma}[Multipliers control for \BALM]\label{thm:penalty:bound}
    Let \cref{assump:prob} hold and suppose that the sequence \(\seq{x^k, \lambda^k, \mu^k}\) is generated by \BALM. Then, for any $k \in \N $:
    \begin{align}
        \tfrac{1}{2\rho_{k+1}}\|\lambda^{k+1}\|^2 + \tfrac{1}{2\nu_{k+1}}\|\mu^{k+1}\|^2 \le c  + \tfrac{1}{2\rho_k}\|\lambda^k\|^2 + \tfrac{1}{2\nu_k}\|\mu^k\|^2,\label{eq:ALM:conv:penalty:bound}
    \end{align}
    where $c\coloneqq f(x^0) - f^\star$. In particular, the following holds for any \(K \in \N\):
    \begin{align*}
        \tfrac{1}{\rho_K}\|\lambda^{K}\|^2 + \tfrac{1}{\nu_K}\|\mu^{K}\|^2 \le 2Kc + \tfrac{1}{\rho_0}\|\lambda^0\|^2 + \tfrac{1}{\nu_0}\|\mu^0\|^2.
    \end{align*}
\end{lemma}
\begin{proof}
    See \appref{app:supp_proofs}.
\end{proof}
Similar to \PBALM, the inequalities \eqref{eq:ALM:L-bound} and \eqref{eq:ALM:conv:penalty:bound} are central to the proof of the successful termination of \BALM. Observe that, in the absence of the proximal term (which is the precise case for \BALM), one gets from \cref{thm:proxALM:conv-step-1} that
\begin{align}
    \|x^{k+1} - \prox_{f_2}(x^{k+1} - \nabla_{x}\calL( x^{k+1},\lambda^{k+1}, \mu^{k+1}))\|_\infty \leq \varepsilon\notag
\end{align}
holds for all $k$ large enough and any $\varepsilon > \tau \in \rnn$, where $\tau \coloneqq \lim_k \tau_k$. The proof follows similarly as the proof of \cref{thm:proxALM:conv-step-1}, with the difference that now the primal residual term \(\tfrac1{\gamk}(x^{k+1} - x^k)\) is absent.
 The same observation is made in the proof of \cref{thm:epsKKT_conv} from which we can establish for \BALM{} that $\|h(x^{k+1})\|^2 \to 0$ and $\|E^{k+1}\|^2 \to 0$ along some infinite subsequence $\calK \subseteq \N$. Bringing these together leads to the following corollary to \cref{thm:ALM:convergence} which states that \BALM{} terminates successfully at an $\varepsilon$-KKT point of problem \eqref{eq:prob}.
\begin{corollary}[$\varepsilon$-KKT termination for \BALM]
    Let \(\seq{x^k, \lambda^k, \mu^k}\) be the sequence of iterates generated by \BALM, and let \cref{assump:prob,assump:penalty} hold. Then, an $\varepsilon$-KKT point of~\eqref{eq:prob} is obtained after finitely many iterations, for any $\varepsilon \in \rp$.
\end{corollary}

\section{Numerical experiments}\label{sec:experiments}
In this section, we present numerical experiments to demonstrate the performance of our proposed schemes. A Python implementation is available online\footnote{\url{https://github.com/adeyemiadeoye/fs-palm}}.
\subsection{Algorithmic setup}\label{sec:setup}
The algorithms considered for comparison are summarized in \cref{tab:algorithms}. Every algorithm is run from the same initial point and with the same inner tolerance sequence $\seq{\tau_k}$.
 
\begin{table}[h]
	\centering
	\caption{Algorithms used for comparison.}
	\label{tab:algorithms}
	\begin{tabularx}{\textwidth}{l X}
		\toprule
		Algorithm & Description \\[1pt]
		\midrule
		
		\PBALM-$\alpha$
		& $\cref{alg:PBALM}$ with $\xi_1=\xi_2=1$ and
		$\phi(k)=(k+1)^\alpha$, where $\alpha>1$.\\[6pt]
		
		\BALM-$\alpha$
		& $\cref{alg:BALM}$ with $\xi_1=\xi_2=1$ and
		$\phi(k)=(k+1)^\alpha$, where $\alpha>1$.\\[6pt]
		
		\texttt{ALM-$\xi$}
		& Standard ALM corresponding to $\xi_1=\xi_2=\xi>1$ and
		$\phi(k)=0$ in $\cref{alg:BALM}$,
		without the initial feasibility condition and the check
		in $\cref{state:ALM:barx}$.\\
		
		\bottomrule
		
		\multicolumn{2}{p{\dimexpr\textwidth-1em\relax}}{
			{\footnotesize \textit{Note.} Unless stated otherwise, we use $\alpha\in\{1.01,2,4\}$ throughout (\cref{sec:setup}). \texttt{ALM-$\xi$}
				takes $\phi\equiv0$ and therefore falls outside \cref{assump:penalty}; it is
				included as a standard-practice baseline and is not covered by our analysis.}
		}
	\end{tabularx}
\end{table}

\subsubsection*{Inner solvers}

Each experiment is run with two inner solvers and is reported for both.
The first is a proximal-gradient method, in which the stepsize is selected by a
standard backtracking linesearch on the smooth part, warm started at $1.1$
times the last accepted stepsize. This inner solver satisfies the hypothesis of
\cref{thm:proxALM:L-bound} by construction. We allow it maximum $5000$ inner
iterations. The second is \texttt{PANOC} \cite{stella2017simple} from the software Alpaqa \cite{pas2022alpaqa} with L-BFGS directions, setting a memory size of $20$ and a maximum iteration number of $1000$. Its linesearch is performed on the forward-backward envelope, so the descent property required by \cref{thm:proxALM:L-bound} is not formally guaranteed at every inner iteration. We did not observe any adverse effect from this in our experiments.

\subsubsection*{Choice of the initial parameters}

In all our tests, we use  $\beta = 0.5$, $\tau_k = \max\{\tau_0\kappa^k,\,\varepsilon\}$ with $\tau_0=0.1$, $\kappa=0.5$ and $\varepsilon$ the requested accuracy, as in \cref{alg:PBALM}, and $\phi_\gamma(k) = (k+1)^{1.01}$ as suggested in \cref{assump:penalty}, unless otherwise stated. For simplicity, we set $\hat{\rho}=\rho_0$, $\hat{\nu}=\nu_0$, and $\hat{\gamma}=\gamma_0$.  The initial multipliers are set to $\lambda^0=\mu^0=0$.

For the initial penalty parameters we adopt the following simple rule of thumb,
\begin{equation}\label{eq:rule3}
  \rho_0 \;=\;
  \frac{\kappa_\rho}{\phi(K_{\mathrm r}+1)}\,
  \frac{\|\nabla^2 f_1(x^0)\|}{\|J_h(x^0)\|^2},
  \qquad
  \nu_0 \;=\;
  \frac{\kappa_\rho}{\phi(K_{\mathrm r}+1)}\,
  \frac{\|\nabla^2 f_1(x^0)\|}{\|J_g(x^0)\|^2},
\end{equation}
with $\kappa_\rho=10^{2}$ and $K_{\mathrm r}=5$. Considering the potential growth according to \(\phi(k+1)\), each penalty modulus reaches
a reasonable reference level after $K_{\mathrm r}$ outer iterations for every $\alpha$.
The dimensionless factor $\kappa_\rho$ is unchanged throughout.
All norms are obtained by power iteration on Hessian- and Jacobian-vector
products, without forming either matrix.
The motivation for this initialization is most direct for $\rho_0$. Since \cref{state:proxALM:barx}
takes $x^0$ feasible, so that $h(x^0)=0$, the term $\tfrac{\rho}{2}\|h(x)\|^2$
contributes exactly $\rho J_h^\top J_h$ to the Hessian of the augmented
Lagrangian at $x^0$, and the corresponding ratio in \eqref{eq:rule3} is the
value of $\rho$ at which the equality constraints are felt as strongly as the
curvature of the objective. We use an expression of the same type for $\nu_0$, with $J_g$
in place of $J_h$. The same reasoning does not carry over, nevertheless, taking the full Jacobian regardless errs towards a smaller initial
penalty, which the growth schedule then recovers.

For the proximal parameters we set $\gamma_0=\kappa_\gamma/\hat L$ with
$\kappa_\gamma=10^{3}$ and $\hat L$ the spectral norm of the Hessian of the
augmented Lagrangian at $x^0$, and $\delta=1/\max\{1,|f_1(x^0)|\}$. Since
$\tfrac{1}{2\gamma}\|x-\hat x\|^2$ adds $1/\gamma$ to the curvature seen by the
subproblem solver, $\kappa_\gamma$ is the inverse fraction of the local
curvature contributed by the proximal term.

\subsubsection*{Initial feasibility and stopping criteria}
\label{sec:stopping-criteria}
The algorithms are stopped when the $\varepsilon$-KKT residual of \cref{def:epsKKT} falls below $\varepsilon=10^{-6}$, that is, when $\max\{\|R_\gamma(x^k)\|_\infty, \|h(x^k)\|_\infty, \|E^k\|_\infty\}\le 10^{-6}$, where $R_\gamma$ denotes the fixed-point residual at unit stepsize. The first term is the primal stationarity measure and is included in the test, so that all three conditions of \cref{def:epsKKT} are certified rather than feasibility and complementarity alone. 
In our tests, we track the total infeasibility measure and the suboptimality gap given, respectively, by $\|h(x^k)\|_\infty + \|\max\{0, g(x^k)\}\|_\infty$ and  $\frac{|f(x^k) - f_{\rm ref}|}{|f(x^0) - f_{\rm ref}|}$, where $f_{\rm ref}$ is a reference value for the optimal value of the problem at hand, specified for each problem below.
For problems where an initial feasible point is not readily available, we solve a phase I problem (see \appref{sec:phase-I}). Note that, to ensure a fair comparison, we use this initial point for all algorithms (including \texttt{ALM}, even though the algorithm does not require a feasible initial point). In our experiments, this is only needed for the CUTEst instances of \Secref{sec:mm_qps}, since for all other problems a feasible point is known by construction. The phase I problem itself always has a trivially feasible point, so \PBALM{} is used to solve it. We refer to \appref{sec:phase-I} for the details and for the limitations of this approach.

In all figures, cost is measured by the cumulative number of evaluations of the gradient of the smooth part of the subproblem objective, including all evaluations performed by the inner solver.

\subsection{Sparse nonconvex QCQP}
\label{sec:qcqp}

We next consider a sparse quadratically constrained quadratic program (QCQP) of the form
\begin{equation}
	\label{eq:qcqp_l1}
	\begin{aligned}
		\minimize_{x\in\rr^n} \quad & \tfrac12 x^\top P_0 x + q_0^\top x + \eta\|x\|_1 \\
		\stt{} \quad & \tfrac12 x^\top P_i x + q_i^\top x - r_i \le 0, \quad i = 1,\dots,m,
	\end{aligned}
\end{equation}
where $P_i\in\rr^{n\times n}$ are symmetric, $q_i\in\rr^n$, and $r_i\in\rr$. No definiteness is assumed on the $P_i$, so the problem is in general nonconvex. The $\ell_1$ penalty parameter is set to $\eta = \tfrac12 \|q_0\|_\infty$.
Problem \eqref{eq:qcqp_l1} is an instance of \eqref{eq:prob} with
\(f_1(x) = \tfrac12 x^\top P_0 x + q_0^\top x\), \(f_2(x) = \eta\|x\|_1\), \(h \equiv 0\), and
\(g(x) = \big(\tfrac12 x^\top P_i x + q_i^\top x - r_i\big)_{i=1}^m \in \rr^m\).

\begin{figure}[t!]
	\centering
	\subfloat{\resizebox*{11cm}{!}{\includegraphics{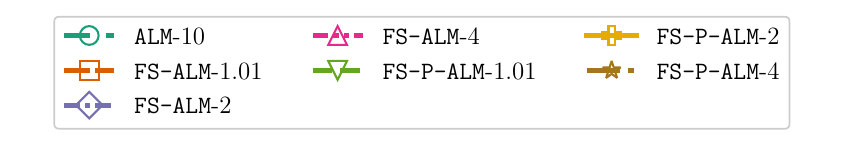}}}
	\vspace{-0.2cm}
	\subfloat[total infeasibility, \texttt{PANOC}]{
		\resizebox*{6.2cm}{!}{\includegraphics{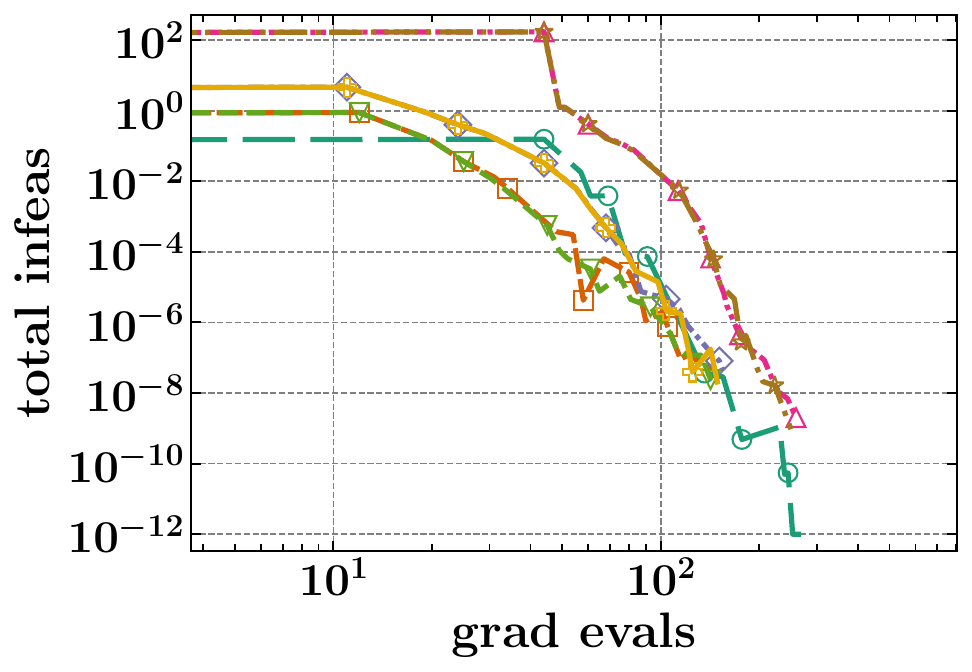}}}\hfil
	\subfloat[objective gap, \texttt{PANOC}]{
		\resizebox*{6.2cm}{!}{\includegraphics{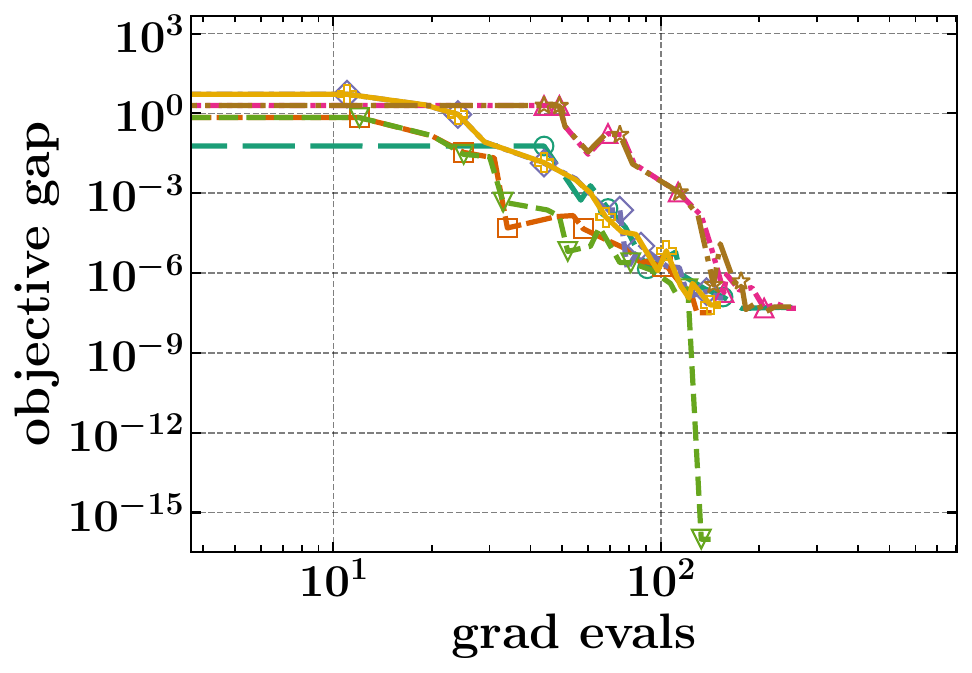}}}
	\vfil
	\subfloat[total infeasibility, \texttt{ProxGrad}]{
		\resizebox*{6.2cm}{!}{\includegraphics{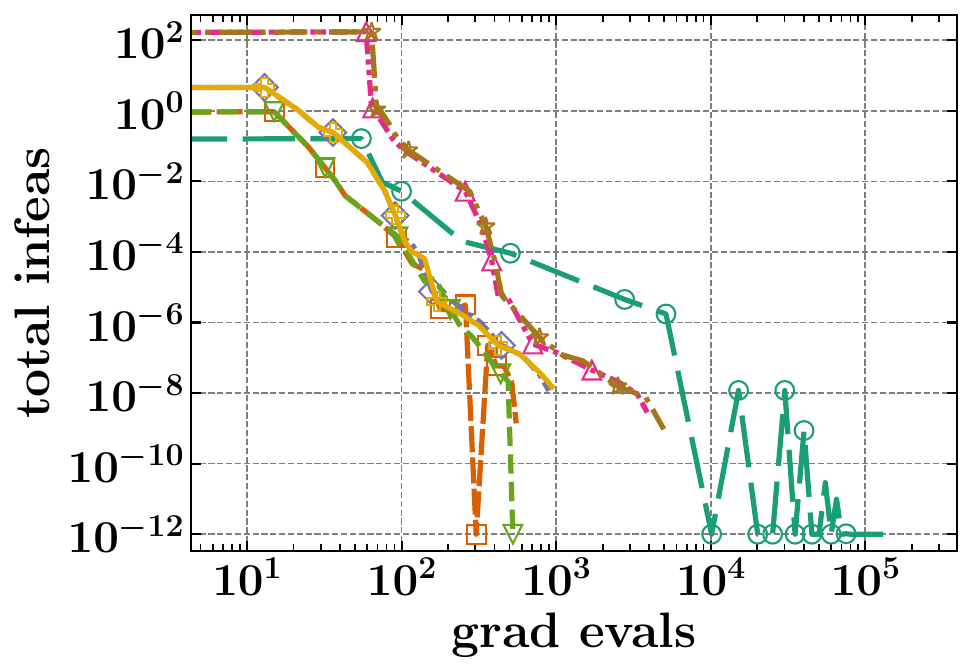}}}\hfil
	\subfloat[objective gap, \texttt{ProxGrad}]{
		\resizebox*{6.2cm}{!}{\includegraphics{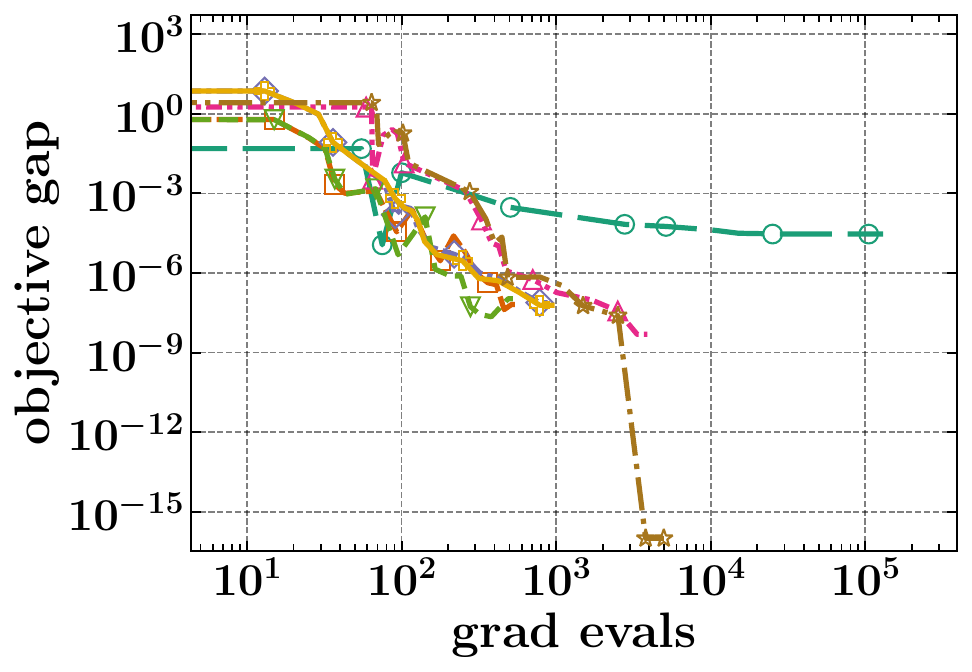}}}
	\caption{Sparse nonconvex QCQP with $n=50$, $m=10$, under both
	inner solvers, for $\alpha\in\{1.01,2,4\}$ (and $\xi=10$ for the \texttt{ALM} variant). Under
	\texttt{PANOC} every run converges. Under \texttt{ProxGrad} \PBALM{} and
	\BALM{} still converge at every $\alpha$, whereas
	\texttt{ALM} does not.}
	\label{fig:qcqp}
\end{figure}

Following the random-instance conventions of \cite{park2017general}, the entries of $q_0,\dots,q_m$ are drawn i.i.d.\ from $\calN(0,1)$. For the objective we form $S = \tfrac{1}{2\sqrt n}(B+B^\top)$ with $B$ having i.i.d.\ $\calN(0,1)$ entries and set $P_0 = V\!\operatorname{diag}(\sigma(w/\|w\|_\infty + \theta))V^\top$, where $S = V\!\operatorname{diag}(w)V^\top$ is the spectral decomposition. The shift $\theta$ and the scale $\sigma$ are drawn per instance so that both the curvature of $f_1$ and the proportion of negative eigenvalues of $P_0$ vary across the test set, $P_0$ remaining indefinite throughout. For the constraints we take $P_i = \sigma_i\,\tfrac1n A_i^\top A_i \succeq 0$ with $A_i\in\rr^{n\times n}$ having i.i.d.\ $\calN(0,1)$ entries and $\sigma_i$ an instance-dependent scale, and $r_i>0$. Each constraint then describes an ellipsoid, their intersection is bounded, and
\(
	x^0 = 0
\)
is strictly feasible for every instance, which supplies the feasible initialization required by \PBALM{} without recourse to the phase~I problem.
Since the problem is nonconvex, its global optimal value is not available. The reference value $f_{\rm ref}$ is therefore the smallest value of $f$ attained at the final points of all runs with the same inner solver in \cref{fig:qcqp} and at a local solution computed by \texttt{IPOPT} \cite{wachter2006implementation}, applied to the standard smooth reformulation of the $\ell_1$ term by means of auxiliary variables.
A sample instance is shown in \cref{fig:qcqp}. Under \texttt{PANOC} the three methods converge comparably, while under \texttt{ProxGrad} only \PBALM{} and \BALM{} reach the tolerance and \texttt{ALM} stalls.

\subsection[Quadratic programs]{Quadratic programs}
\label{sec:mm_qps}
We illustrate the proposed algorithms on quadratic programming (QP) problems  in the form \eqref{eq:prob} with $f_1(x) = \frac{1}{2} x^\top Q x + q^\top x + c$, where $Q\in\rr^{n\times n}$ is symmetric and not necessarily positive semidefinite, $q\in\rr^n$, $c\in\rr$; $f_2(x) = \delta_{\Omega}(x)$, where $\Omega = \{x\in\rr^n: x_l \le x \le x_u\}$; and $g(x) = [Ax - u, l - Ax]$, where $A\in\rr^{m\times n}$.

We report two sets of experiments, one on instances taken from an
established public collection and one on randomly generated instances whose
conditioning we control.

\subsubsection*{CUTEst instances}

\begin{figure}[t!]
	\centering
	\resizebox*{\linewidth}{!}{\includegraphics{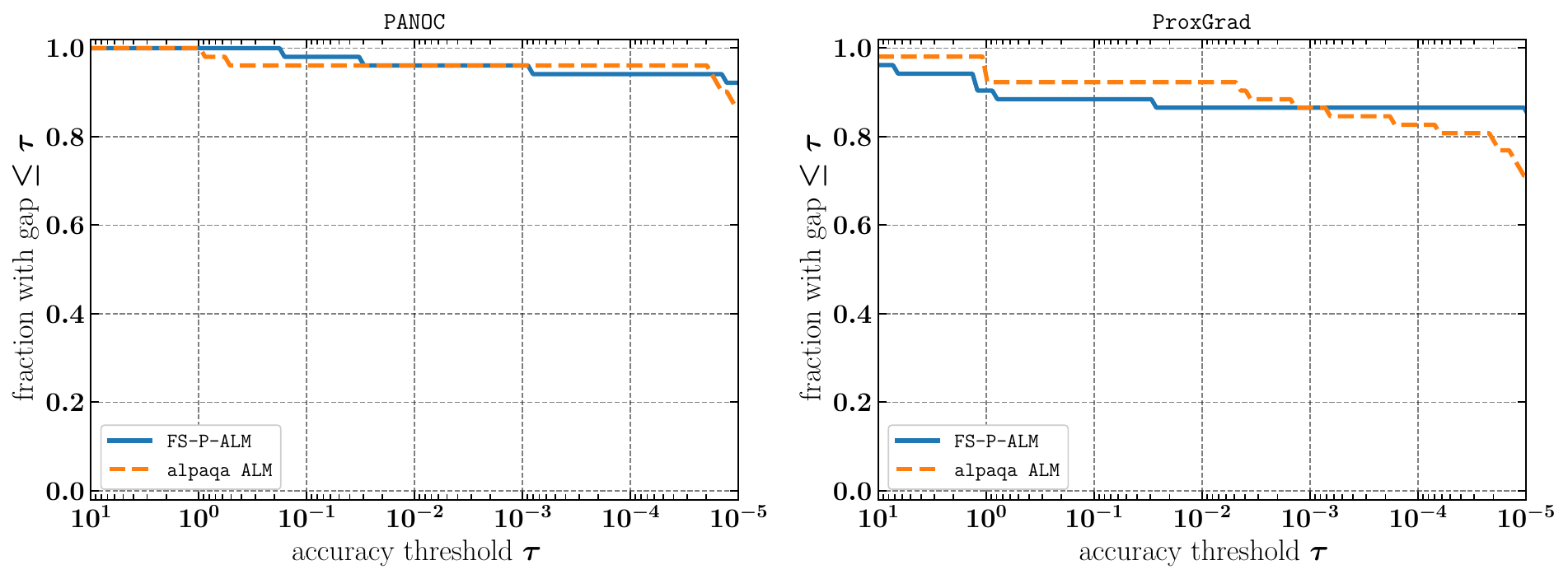}}
	\caption{Accuracy attained by the point each method returns, against the
	threshold $\tau$, on the $53$ CUTEst QPs.
Under \texttt{PANOC}, \PBALM{} and \texttt{alpaqa} satisfy their
respective stopping tests on $46$ and $47$ instances; under
\texttt{ProxGrad}, on $36$ and $34$. 
	}
	\label{fig:cutest-profile}
\end{figure}

The first set is drawn from the CUTEst
collection \cite{gould2015cutest} and comprises $53$ instances. Of these, $44$
are convex, from the Maros-M{\'e}sz{\'a}ros repository
\cite{maros1999repository} and related CUTEst QP families, while the remaining
nine are the \texttt{NCVXQP} instances, whose Hessians are indefinite. The dimensions range from $n=2$ to $n=745$ and from
$m=1$ to $m=1001$. The upper end of this range is a limitation of our 
implementation rather than of the method itself, since it assembles Hessians
and Jacobians densely for automatic differentiation and therefore cannot reach
the largest sparse instances of the collection. 

On this set we compare \PBALM{} against the augmented Lagrangian solver of
\texttt{alpaqa} \cite{pas2022alpaqa}, which safeguards the multipliers by
projection onto a bounded set and increases the penalty multiplicatively. Both
methods drive the same \texttt{PANOC} inner solver with identical settings and
run in the same process. We use $\alpha=4$ and $\alpha_\gamma=4$, a common stopping tolerance of
$10^{-5}$ for both methods, and a limit of $500$ outer iterations. The tolerance
enters only through where each method chooses to stop, and hence which point it
returns.

\Cref{fig:cutest-profile} reports, against the accuracy threshold $\tau$, the
fraction of instances on which the point a method \emph{returns} attains
$\max\{\text{violation},\,|f-f_{\rm ref}|/\max\{1,|f_{\rm ref}|\}\}\leq\tau$,
where $f_{\rm ref}$ is the best value either method reached at a feasible point. 
In this figure we report the attained accuracy rather than the cost of reaching a prescribed target because the two methods do not stop on the same test. 
Neither method dominates over the whole range, and the ordering depends both on
the threshold and on the inner solver. 
Under \texttt{PANOC}, \PBALM{} satisfies its stopping test on $46$
instances and \texttt{alpaqa} on $47$; under \texttt{ProxGrad}, the
corresponding counts are $36$ and $34$. Among the remaining runs,
\PBALM{} reaches the outer-iteration limit while returning a point
feasible to $10^{-5}$ on $4$ instances under \texttt{PANOC} and $11$
under \texttt{ProxGrad}. The corresponding counts for \texttt{alpaqa}
are $5$ and $17$. These runs are feasible according to the common post-processing measure, but do not satisfy the method's complete internal stopping test and hence are not certified as converged.
\PBALM{} returns a genuinely infeasible point on $1$ instance under
\texttt{PANOC} and $5$ under \texttt{ProxGrad}; the corresponding
counts for \texttt{alpaqa} are $1$ and $2$. On \texttt{QPCBOEI1}, and
under \texttt{ProxGrad} also on \texttt{QPCBOEI2}, phase~I does not
produce a feasible starting point, so \PBALM{} is not run.

\begin{figure}[t!]
	\centering
	\resizebox*{\linewidth}{!}{\includegraphics{ 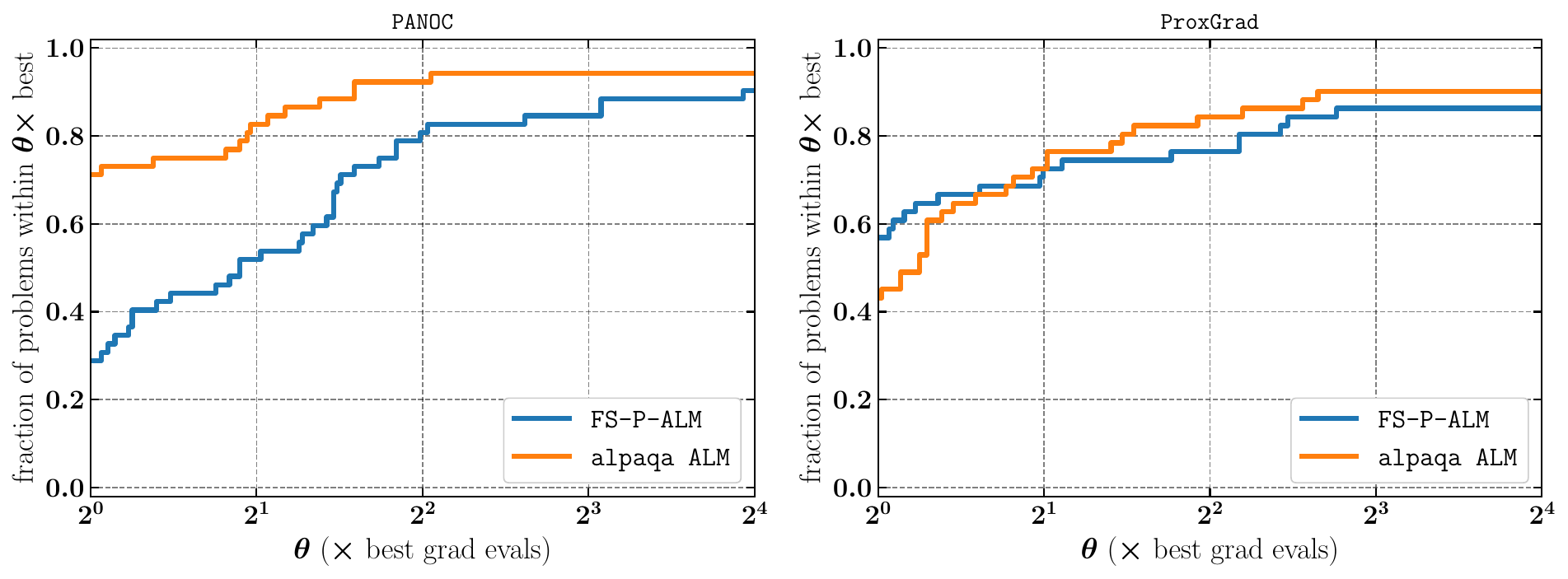}}
	\caption{Performance profile in gradient evaluations on the same set,
	in the sense of \cite{dolan2002benchmarking}. 
	The ordering reverses with the inner solver. Under
	\texttt{PANOC}, whose quasi-Newton directions absorb subproblem
	ill-conditioning, \texttt{alpaqa} is cheaper on $35$ of the $48$ instances
	both solve. Under \texttt{ProxGrad}, where that conditioning is felt
	directly, \PBALM{} is cheaper on $26$ of $42$.}
	\label{fig:cutest-cost}
\end{figure}

	\Cref{fig:cutest-cost} plots the performance profile \cite{dolan2002benchmarking} in gradient evaluations on the same dataset. 
	Cost is the number of
	evaluations at which a method first reaches one common target, feasible to
	$10^{-5}$ and within $10\%$ of the best objective either method found.
	Thus, a higher curve means that the
method reaches the common target within a smaller multiple of the best
recorded cost on more jointly successful problems. Under \texttt{PANOC},
\texttt{alpaqa} has the smaller cost on $35$ of the $48$ jointly successful
instances; under \texttt{ProxGrad}, \PBALM{} has the smaller cost on $26$ of
$42$. 
This confirms the overall observation that Panoc directions absorb the ill-conditioning of the subproblems, while ProxGrad does not.

\subsubsection*{Randomly generated instances}

\begin{figure}[t!]
	\centering
	\subfloat{\resizebox*{11cm}{!}{\includegraphics{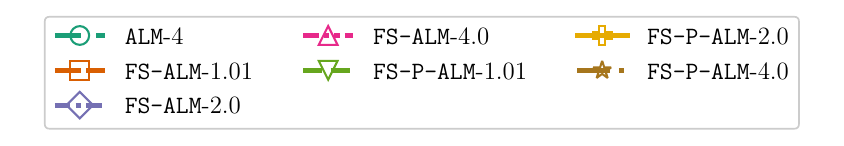}}}
	\vspace{-0.2cm}
	\subfloat[total infeasibility, \texttt{PANOC}]{
		\resizebox*{6.2cm}{!}{\includegraphics{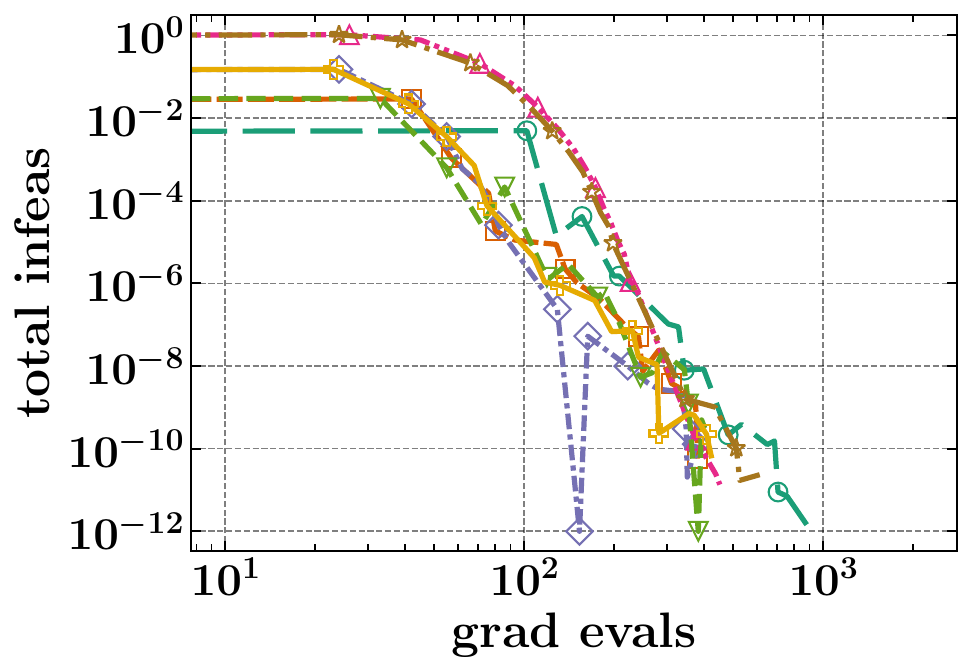}}}\hfil
	\subfloat[objective gap, \texttt{PANOC}]{
		\resizebox*{6.2cm}{!}{\includegraphics{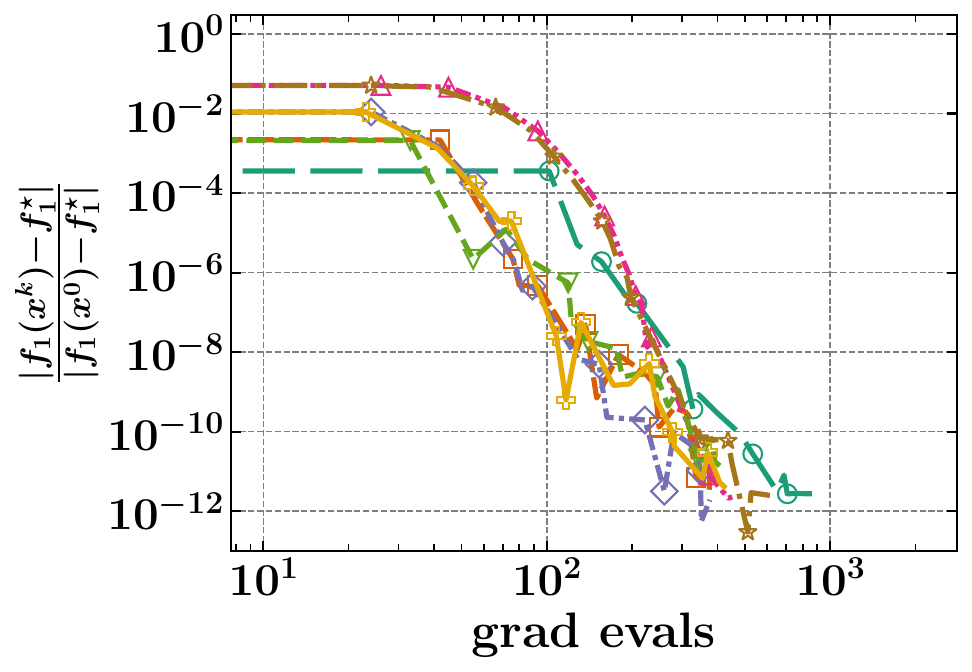}}}
	\vfil
	\subfloat[total infeasibility, \texttt{ProxGrad}]{
		\resizebox*{6.2cm}{!}{\includegraphics{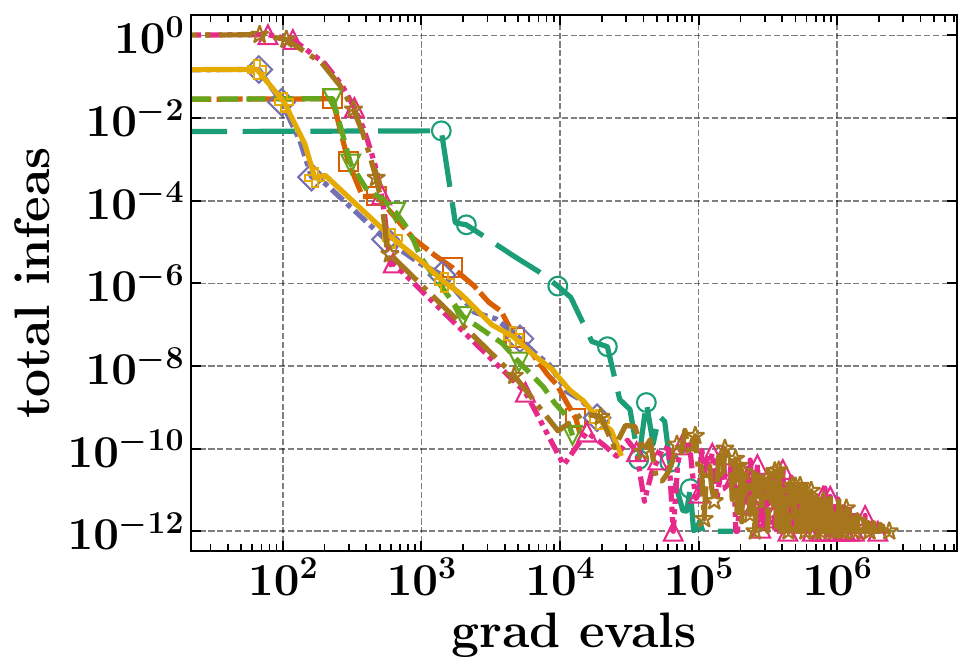}}}\hfil
	\subfloat[objective gap, \texttt{ProxGrad}]{
		\resizebox*{6.2cm}{!}{\includegraphics{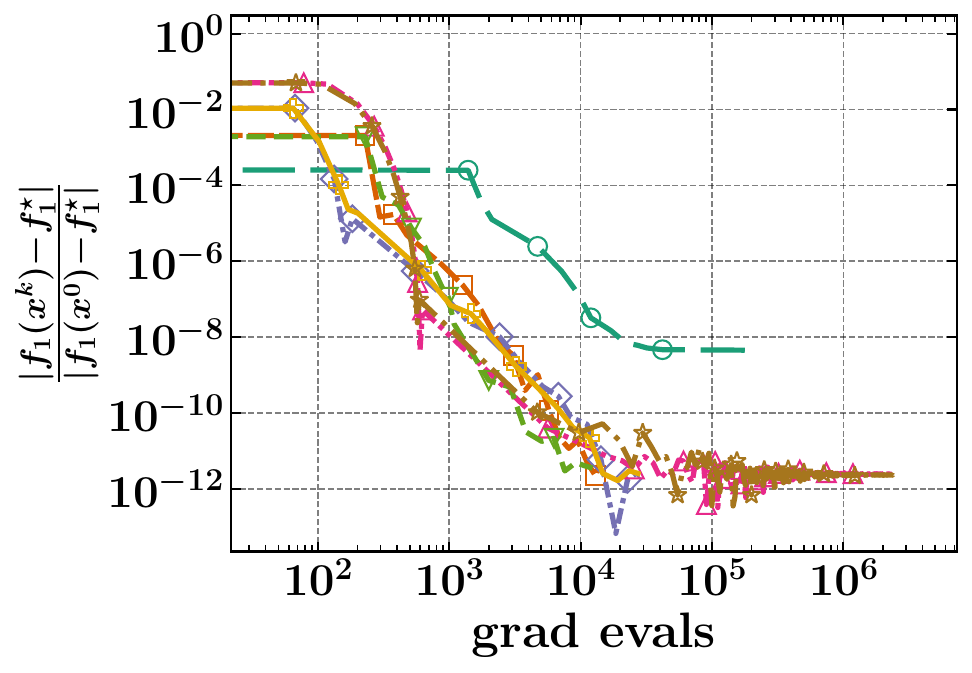}}}
	\caption{Convex QP with $n=60$, $m=30$ and $\kappa(Q)=10^{2}$, under both
	inner solvers, for $\alpha\in\{1.01,2,4\}$ (and $\xi=4$ for the \texttt{ALM} variant). Under
	\texttt{PANOC} every run converges. Under \texttt{ProxGrad}, where the
	subproblem conditioning is felt directly, \texttt{ALM} breaks down while
	\PBALM{} and \BALM{} still converge.}
	\label{fig:ALM-comparison-00}
\end{figure}

The second set consists of randomly generated convex QPs, which allows
the conditioning of the subproblems to be varied in a controlled way. The
instances are generated with a prescribed condition number $\kappa(Q)$. Concretely, $Q = V\diag(w)V^\top$ with $V$ a random orthogonal matrix and $w$ log-spaced over $[1,\kappa(Q)]$. $q$ has i.i.d.\ standard Gaussian entries and $c=0$. The linear constraints $A$ are placed around a strictly feasible interior point, with a quarter of the rows tightened to near-activity there so the constraints are actually exercised rather than redundant, and the box is fixed at $[-10,10]^n$. The same interior point serves as $x^0$, so no phase I step is needed. Since these instances are convex, the reference value $f_{\rm ref}$ is their optimal value, which we compute with the interior-point solver \texttt{Clarabel} \cite{goulart2024clarabel} through \texttt{CVXPY} \cite{diamond2016cvxpy}. The results are displayed in \cref{fig:ALM-comparison-00}, for $\alpha\in\{1.01,2,4\}$ in accordance with \cref{tab:algorithms}, with the initial parameters set as in \cref{sec:setup}. Under \texttt{PANOC} all runs converge. Under \texttt{ProxGrad}, where the subproblem conditioning is felt directly, \texttt{ALM} does not reach the tolerance while \PBALM{} and \BALM{} do.

\subsection{Basis pursuit (nonconvex formulation)}
\label{sec:basis_pursuit}
We consider the basis pursuit problem which involves finding the sparsest solution to an undetermined linear system of equations. That is,
\begin{equation}
	\label{eq:basis_pursuit1}
\begin{aligned}
    \minimize_{z\in\rr^n} \quad & \|z\|_1 \qquad \stt{} \quad & Bz = b,
\end{aligned}
\end{equation}
where \(B\in\rr^{p\times n}\), with \(p<n\), and \(b\in\rr^p\). We consider a nonconvex reformulation of \eqref{eq:basis_pursuit1} \cite{sahin2019inexact}:
\begin{equation*}
\begin{aligned}
    \minimize_{x\in\rr^{2n}} \quad & \|x\|^2 \qquad\stt{} \quad & \bar{B}x^{\circ 2} = b,
\end{aligned}
\end{equation*}
where $x \coloneqq (x_1, x_2)\in \R^{2n}$, $\bar{B} \coloneqq [B, -B] \in \rr^{p\times 2n}$, and \(\circ\) is used indicate element-wise power. Letting $z^+, z^-$ denote the positive and negative parts of \(z\) satisfying $z = z^+ - z^-$, it holds that 
$x_1^{\circ 2} = z^+$, $x_2^{\circ 2} = z^-$. The entries of $B$ are independent draws from the standard normal distribution $\calN(0,1)$. We then set $b=Bz^\star$, where $z^\star$ is also a random vector with i.i.d. entries from $\calN(0,1)$ having $k$ nonzero values fixed at $10$. The proximal parameters are set as in \cref{sec:setup}.
The reference value is $f_{\rm ref}=\|z^\star\|_1$, which is the optimal value of both \eqref{eq:basis_pursuit1} and its reformulation. That $z^\star$ solves \eqref{eq:basis_pursuit1} is verified by solving its linear programming reformulation with \texttt{HiGHS} \cite{huangfu2018parallelizing}.

\begin{figure}[t!]
	\centering
	\subfloat{\resizebox*{11cm}{!}{\includegraphics{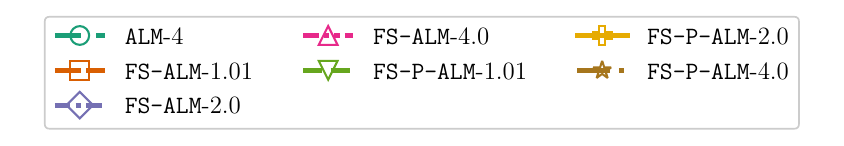}}}
	\vspace{-0.2cm}
	\subfloat[total infeasibility, \texttt{PANOC}]{
		\resizebox*{6.2cm}{!}{\includegraphics{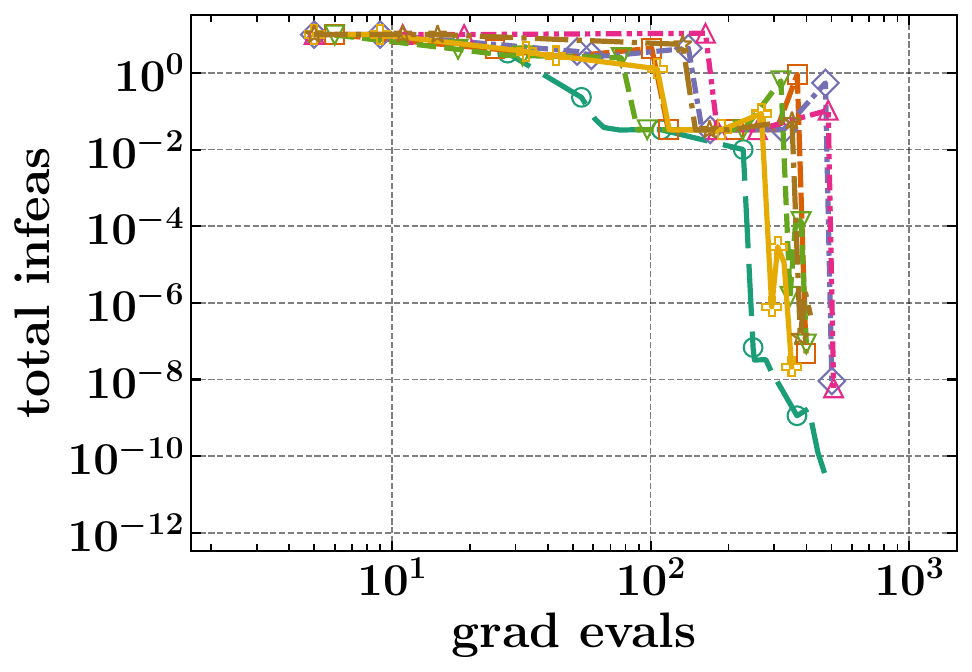}}}\hfil
	\subfloat[objective gap, \texttt{PANOC}]{
		\resizebox*{6.2cm}{!}{\includegraphics{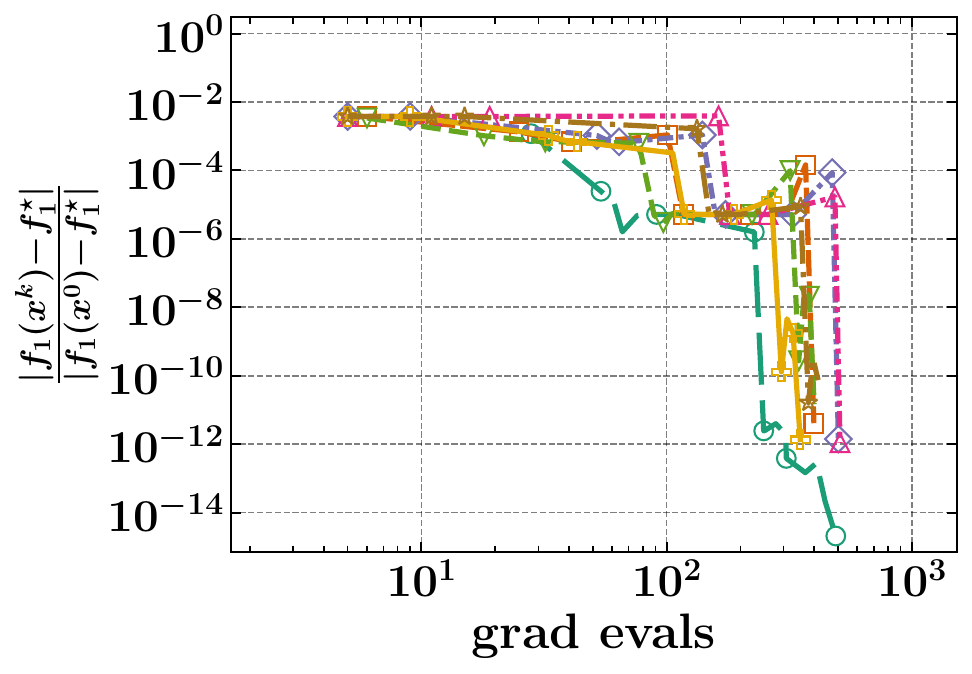}}}
	\vfil
	\subfloat[total infeasibility, \texttt{ProxGrad}]{
		\resizebox*{6.2cm}{!}{\includegraphics{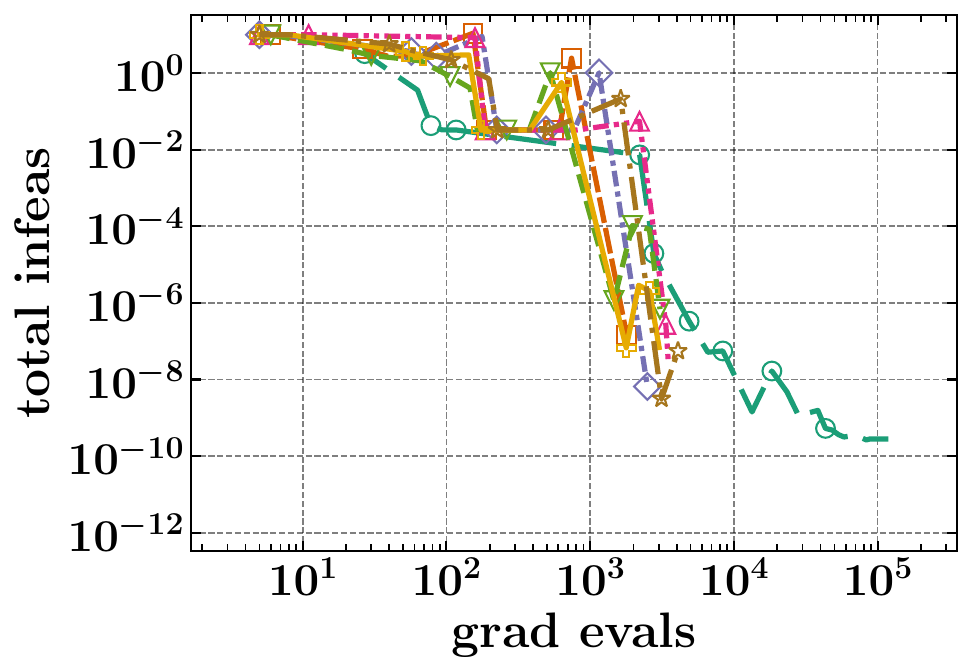}}}\hfil
	\subfloat[objective gap, \texttt{ProxGrad}]{
		\resizebox*{6.2cm}{!}{\includegraphics{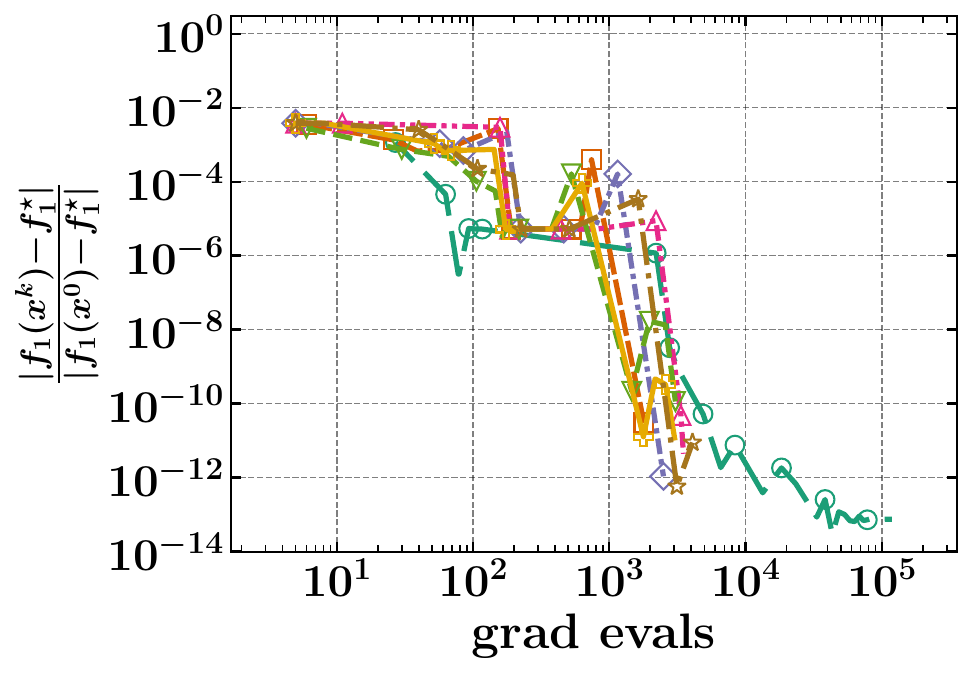}}}
	\caption{Basis pursuit with $p=400$, $n=1024$ and $k=10$, under both inner
	solvers, for $\alpha\in\{1.01,2,4\}$ and $\xi=4$. Under \texttt{PANOC} all
	runs converge and lie within a small factor of one another. Under
	\texttt{ProxGrad} \PBALM{} and \BALM{} still converge at every
	$\alpha$, whereas \texttt{ALM} does not.}
	\label{fig:basis-pursuit}
\end{figure}

\Cref{fig:basis-pursuit} compares the algorithms for $\alpha\in\{1.01,2,4\}$ under both inner solvers. Under \texttt{PANOC} the three methods lie within a small factor of one another, while  under \texttt{ProxGrad} \PBALM{} and \BALM{} converge at every $\alpha$ while \texttt{ALM} does not.
Both \PBALM{} and \BALM{} provide advantages over \texttt{ALM} on these problems, avoiding the ill-conditioning issue. The size of the advantage depends on how the subproblems are solved and on how the initial penalty is chosen. Under a quasi-Newton inner solver, and with each method given a well-scaled initial penalty, the three lie within a small factor of one another. The separation appears when the subproblems are solved by a first-order method, or when the initial penalty is misspecified. See also \appref{sec:supporting} for the effect of \(\alpha\) given a fixed initialization of the penalty parameter.

\begin{figure}[t!]\centering
  \resizebox*{13cm}{!}{\includegraphics{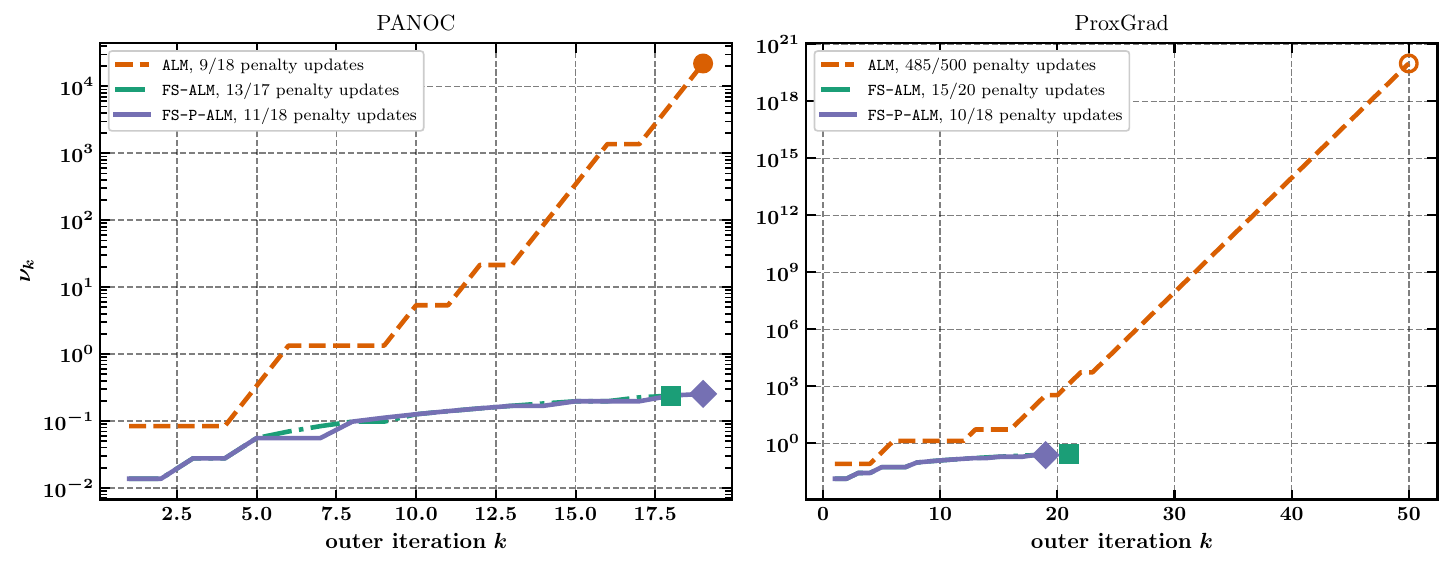}}
  \caption{Penalty $\rho_k$ against outer iteration on a basis-pursuit
  instance, under each inner solver, labelled by the number of $\beta$-test
  firings. Filled end markers denote runs that converged, open markers runs
  that did not. Curves are drawn over the range in which the trajectory
  remains numerically meaningful.}
  \label{fig:aux-runaway}
\end{figure}

\Cref{fig:aux-runaway} shows the penalty trajectories on a basis-pursuit
instance, labelled by the number of times the $\beta$-progress test fires.
That test, rather than an inexact subproblem solve, is what triggers a penalty
update, and it fires at comparable rates for all three methods. \texttt{ALM} multiplies by $\xi$ on each firing, so firings
compound geometrically, whereas our update reduces to
$\rho\leftarrow\max\{\rho,\hat\rho\phi(k+1)\}$ and is bounded by the schedule.

These results demonstrate the promising practical performance of our proposed approach on nonconvex problems, highlighting their potential to effectively solve challenging problems in many application domains.
\section{Conclusions}
\label{sec:conclusion}

We proposed an inexact proximal augmented Lagrangian method that enables aggressive and adaptive penalty updates for solving nonconvex structured optimization problems. Our theoretical analysis is based on properties of the iterates that stem from the initial feasibility assumption, establishing that the method attains an approximate Karush--Kuhn--Tucker point, up to a prescribed accuracy, in finitely many iterations. The analysis further extends to the limiting case without the proximal term, thereby encompassing the classical augmented Lagrangian method as a special case. Numerical experiments on both convex and nonconvex problems of varying complexity highlight the efficiency of the proposed approach.

Future research directions include extending the method to incorporate positive definite matrices in place of scalar proximal stepsizes, together with appropriate update rules for such a matrix. Another important avenue is the development of iteration complexity estimates for the proposed framework.

\section*{Acknowledgments}
The authors thank Pieter Pas (KU Leuven) for informed discussions on the PANOC solver.
\section*{Declarations}
\subsection*{Conflict of interest} The authors have no conflict of interest to declare.
\subsection*{Data/code availability} Some numerical experiments were carried out on a subset of the Maros-M{\'e}sz{\'a}ros quadratic programming test set, which is publicly available at \url{https://www.doc.ic.ac.uk/~im/} and \url{https://old.sztaki.hu/~meszaros/public_ftp/qpdata/} (see also \cite{maros1999repository}). A Python implementation of our approach, together with all supporting scripts, is publicly available on GitHub at \url{https://github.com/adeyemiadeoye/fs-palm}.

\clearpage
\appendix
\section{The phase I problem for initialization}\label{sec:phase-I}
When an initial \emph{feasible} point is not readily available, a so-called phase I problem can be solved to find one for \PBALM. We consider the following problem:
\begin{equation}\label{eq:phase-I-2}
    \begin{aligned}
        \minimize_{(x,s) \in \R^{n+1}} \quad & \tfrac{1}{2}\|h(x)\|^2 + s^2 + \indicator_{\dom f_2}(x) \quad 
        \stt{} \quad &g(x) \leq s\mathbf 1.
    \end{aligned}
\end{equation}
Here, $\mathbf 1\in\R^m$ is the vector of all ones, and $\dom f_2$ is assumed to be closed. For any $x^0\in\dom f_2$, the pair $(x^0,s^0)$ with $s^0=\max\{g_1(x^0),\dots,g_m(x^0)\}$ is feasible for \eqref{eq:phase-I-2}, so that \PBALM{} can be applied to it starting from this point. Moreover, \eqref{eq:prob} is feasible if and only if \eqref{eq:phase-I-2} admits a point $(x,s)$ with zero objective value, in which case $s=0$ and $x$ is feasible for \eqref{eq:prob}. However, since \eqref{eq:phase-I-2} is nonconvex in general, \PBALM{} is only guaranteed to return an $\varepsilon$-KKT point of it, which need not be feasible for \eqref{eq:prob}. The outcome of phase I thus has to be checked. In our experiments, a feasible initial point is available by construction for the QCQP of \Secref{sec:qcqp}, the randomly generated quadratic programs of \Secref{sec:mm_qps}, and the basis pursuit problem of \Secref{sec:basis_pursuit}. Phase I is only needed for the CUTEst instances of \Secref{sec:mm_qps}, where it does not produce a feasible point for \texttt{QPCBOEI1}, and under \texttt{ProxGrad} also for \texttt{QPCBOEI2}, so that \PBALM{} is not run on these instances.
\section{Proofs for \cref{sec:convergence,sec:ALM-algorithm}}\label{app:supp_proofs}
\begin{appendixproof}{thm:E-KKT}
By the definition of $E^k$ in~\eqref{eq:E} and $\|E^{k+1}\|_\infty \leq \varepsilon$, we have that if $E_i^{k+1} = \tfrac{1}{\nu_k}\mu_i^k$, then $g_i(x^{k+1}) \le 0 \le \varepsilon$, since $\tfrac{1}{\nu_k}\mu_i^k \in \rnn$. Likewise, $|E_i^{k+1}| \le \varepsilon$ implies that $g_i(x^{k+1}) \le \varepsilon$ whenever $E_i^{k+1} = -g_i(x^{k+1})$. Hence, $g_i(x^{k+1}) \le \varepsilon$ for all $i \in [m]$ and $\|\max\{0,g(x^{k+1})\}\|_\infty \le \varepsilon$. Now, if $g_i(x^{k+1}) < -\varepsilon$, then by $|E_i^{k+1}| \le \varepsilon$, we must have $E_i^{k+1} = \tfrac{1}{\nu_k}\mu_i^k \le \varepsilon$. These give $g_i(x^{k+1}) + \tfrac{1}{\nu_k}\mu_i^k < - \varepsilon + \varepsilon = 0$. By the definition of $\mu_i^{k+1}$ in \cref{state:proxALM:mu}, we have $\mu_i^{k+1} = 0$.
\end{appendixproof}
\begin{appendixproof}{thm:penalty:bound}
    Using \cref{state:ALM:lambda,state:ALM:mu} of \BALM{} in the definition \eqref{eq:L'} of $\calL_{\rho,\nu}$, we get
    \begin{align*}
        \tfrac{1}{2\rho_k}\|\lambda^{k+1}\|^2 + \tfrac{1}{2\nu_k}\|\mu^{k+1}\|^2 &= \tfrac{\rho_k}{2}\|h(x^{k+1}) + \tfrac{\lambda^k}{\rho_k}\|^2 + \tfrac{1}{2\nu_k}\|\max\{0, \nu_k g(x^{k+1}) + \mu^k\}\|^2\\
        &= \calL_{\rho_{k},\nu_{k}}(x^{k+1}, \lambda^k, \mu^k) - f_1(x^{k+1}) + \tfrac{1}{2\rho_k}\|\lambda^k\|^2 + \tfrac{1}{2\nu_k}\|\mu^k\|^2\\
        &\le f(x^0) - f_2(x^{k+1}) - f_1(x^{k+1}) + \tfrac{1}{2\rho_k}\|\lambda^k\|^2 + \tfrac{1}{2\nu_k}\|\mu^k\|^2\\
        &\le f(x^0) - f^\star  + \tfrac{1}{2\rho_k}\|\lambda^k\|^2 + \tfrac{1}{2\nu_k}\|\mu^k\|^2,
    \end{align*}
    where the first inequality follows from \eqref{eq:ALM:L-bound} and the second one from \cref{assump:obj}. Since $\rho_{k+1}\ge\rho_k$ and $\nu_{k+1}\ge\nu_k$, the left-hand side of \eqref{eq:ALM:conv:penalty:bound} is bounded above by the left-hand side of the above chain, which proves \eqref{eq:ALM:conv:penalty:bound}. Summing \eqref{eq:ALM:conv:penalty:bound} over $k=0,\dots,K-1$ and multiplying by $2$ yields the last claim.
\end{appendixproof}
\section{Additional numerical results}\label{sec:penalty-param-supp-add}\label{sec:supporting}

\subsection*{Effect of the proximal term}

\begin{figure}[h]\centering
  \resizebox*{9cm}{!}{\includegraphics{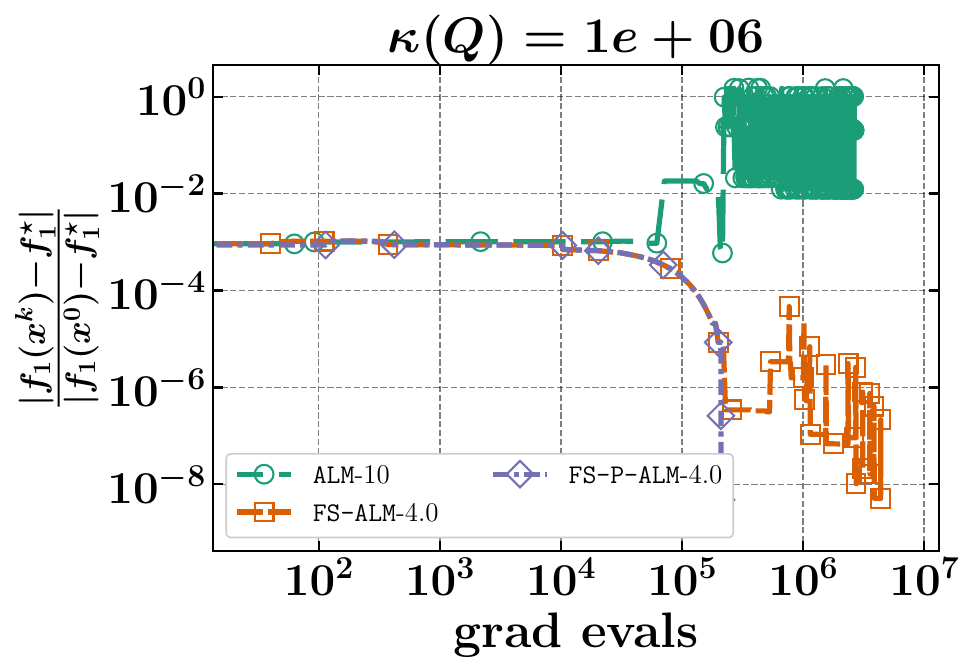}}
  \caption{Relative suboptimality against gradient evaluations on
  \texttt{DUALC1} under \texttt{PANOC}, for \texttt{ALM}, \BALM{} and \PBALM{}
  at a common growth exponent. Only \PBALM{} reaches the requested tolerance.}
  \label{fig:dualc1-prox}
\end{figure}

\begin{figure}[p]
	\centering
	\resizebox*{0.9\linewidth}{!}{\includegraphics{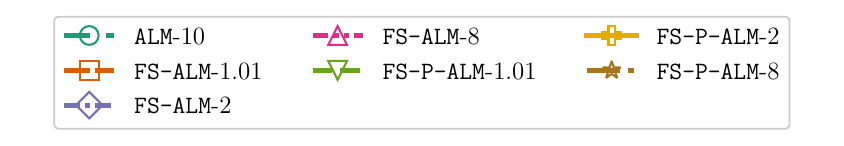}}
	\par\medskip
	\subfloat[CVXQP2\_M, \texttt{PANOC}]{
		\resizebox*{5.6cm}{!}{\includegraphics{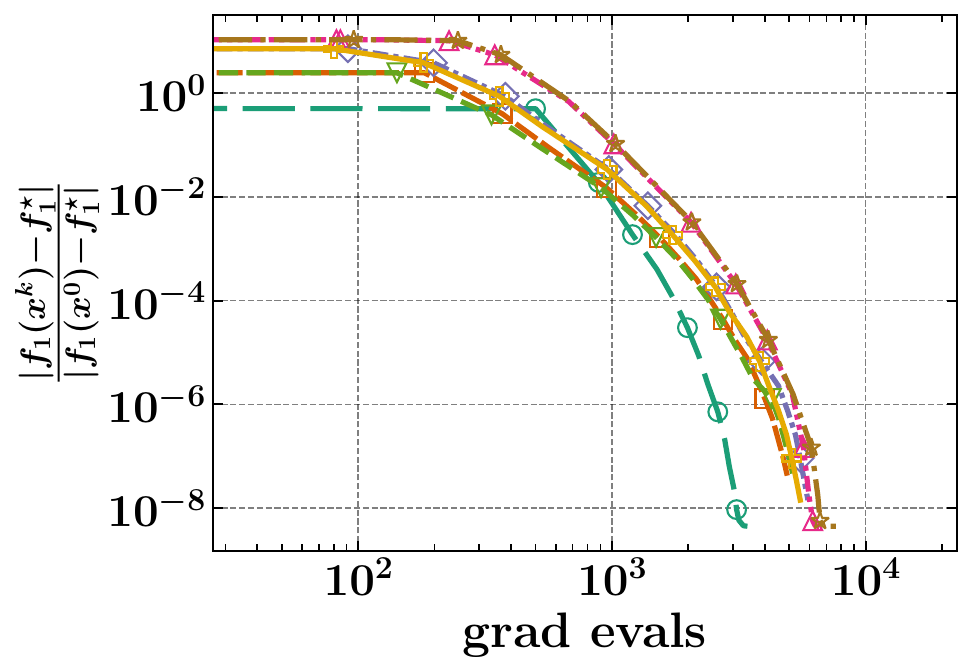}}}\hspace{0.6cm}
	\subfloat[CVXQP2\_M, \texttt{ProxGrad}]{
		\resizebox*{5.6cm}{!}{\includegraphics{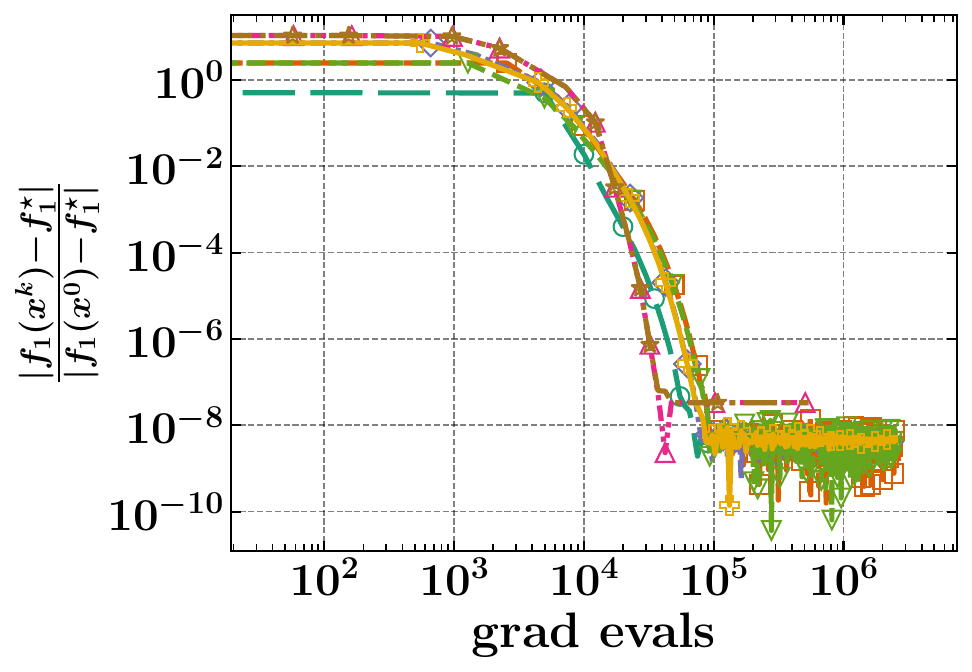}}}
	\vfill
	\subfloat[DUALC1, \texttt{PANOC}]{
		\resizebox*{5.6cm}{!}{\includegraphics{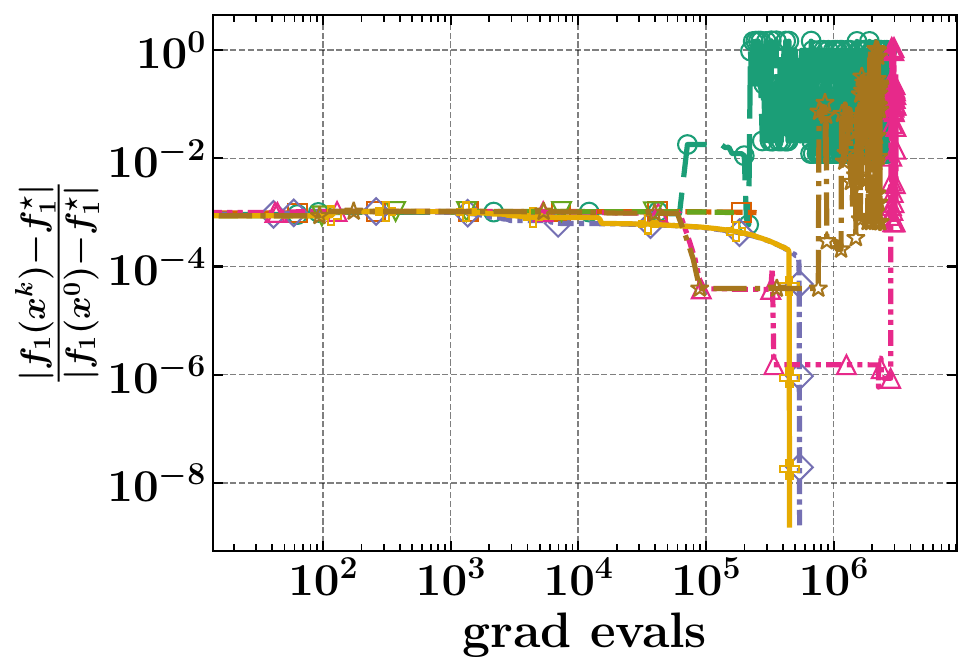}}}\hspace{0.6cm}
	\subfloat[DUALC1, \texttt{ProxGrad}]{
		\resizebox*{5.6cm}{!}{\includegraphics{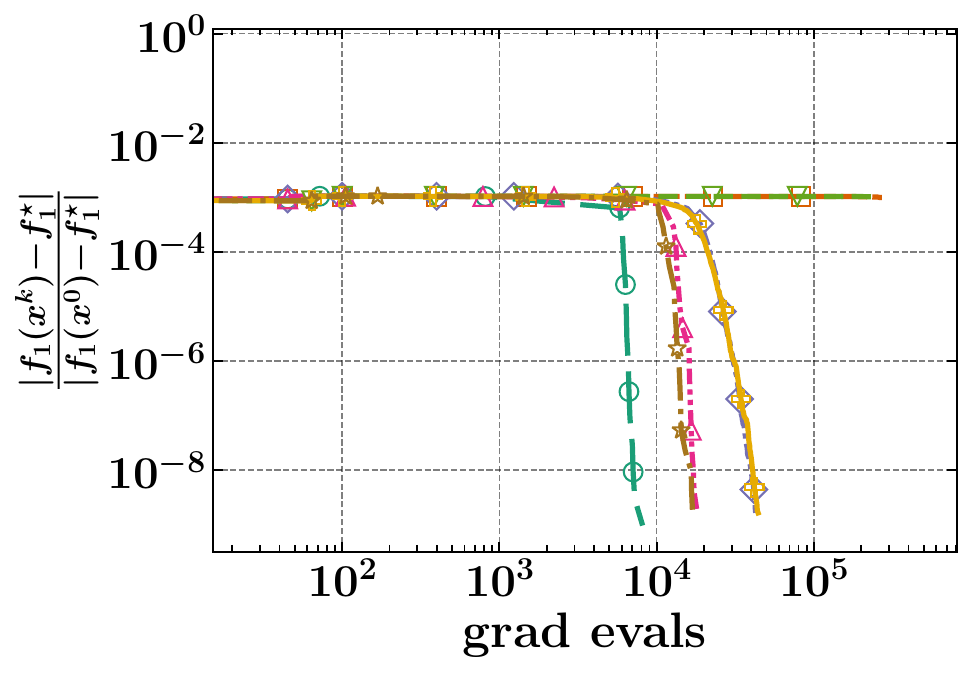}}}
	\vfill
	\subfloat[GENHS28, \texttt{PANOC}]{
		\resizebox*{5.6cm}{!}{\includegraphics{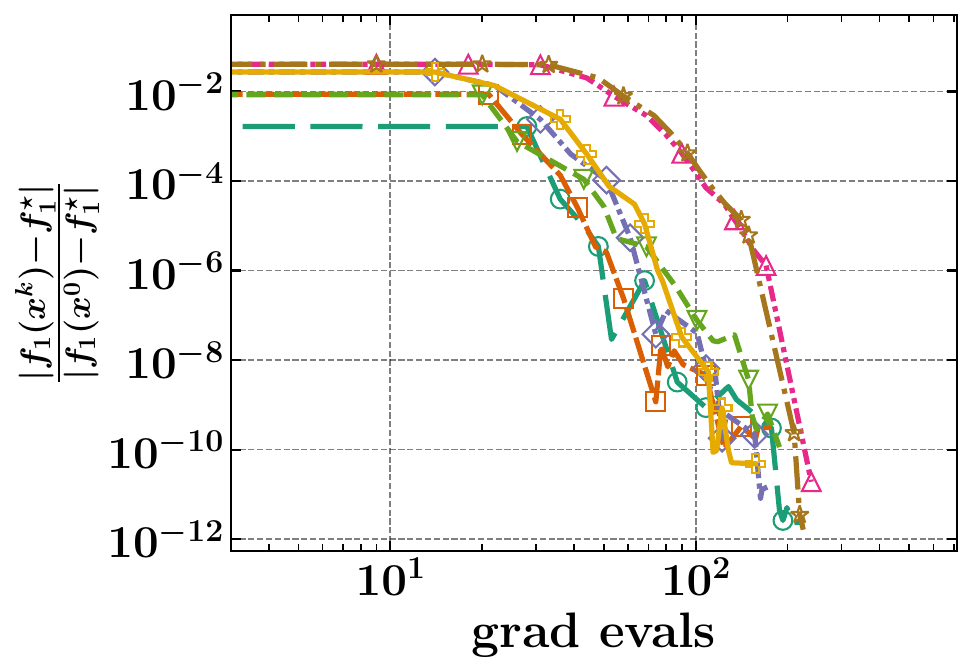}}}\hspace{0.6cm}
	\subfloat[GENHS28, \texttt{ProxGrad}]{
		\resizebox*{5.6cm}{!}{\includegraphics{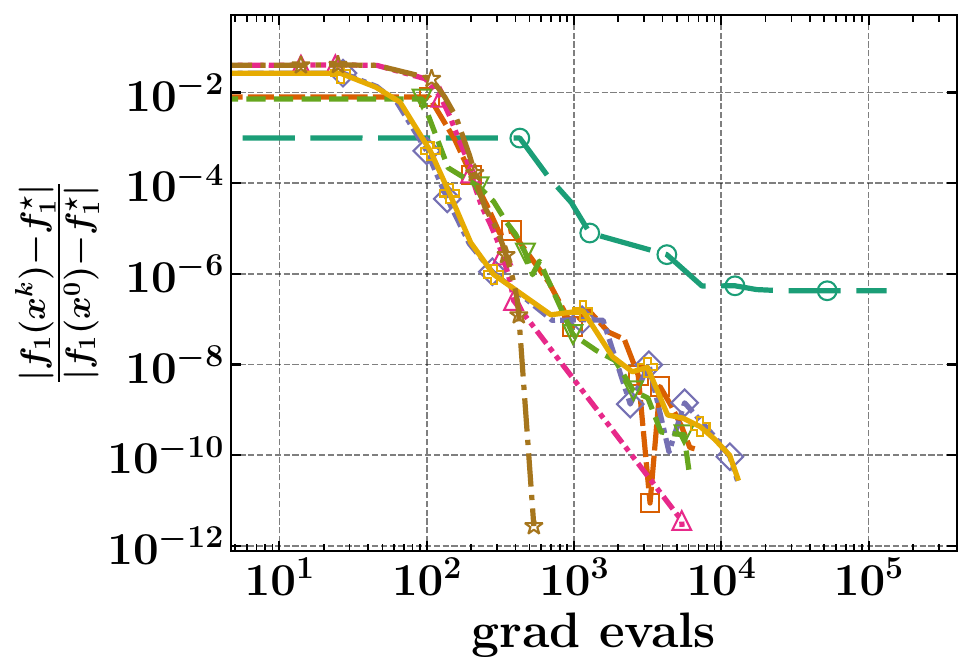}}}
	\vfill
	\subfloat[LOTSCHD, \texttt{PANOC}]{
		\resizebox*{5.6cm}{!}{\includegraphics{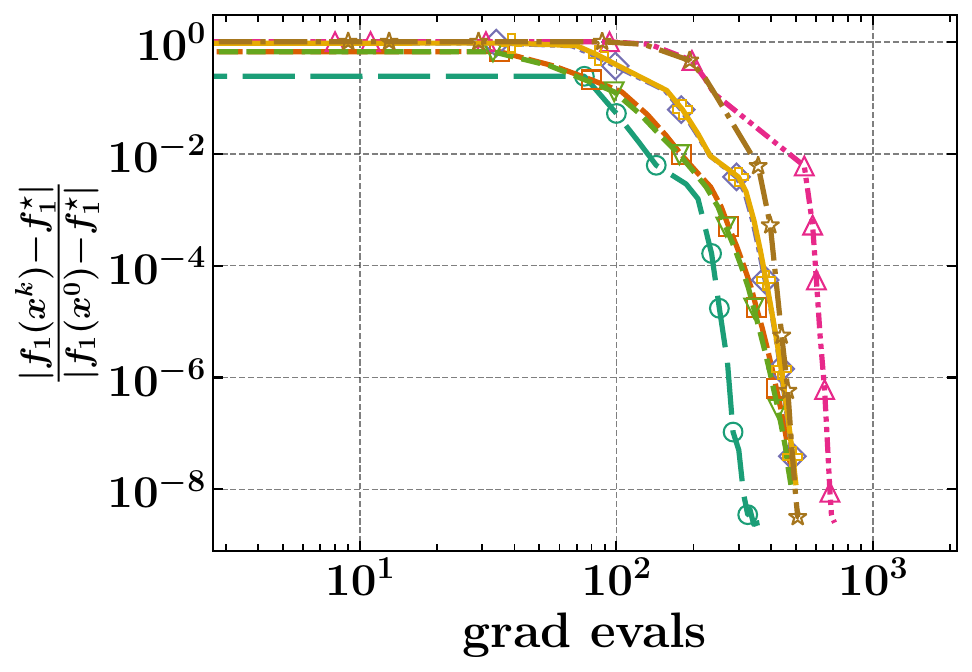}}}\hspace{0.6cm}
	\subfloat[LOTSCHD, \texttt{ProxGrad}]{
		\resizebox*{5.6cm}{!}{\includegraphics{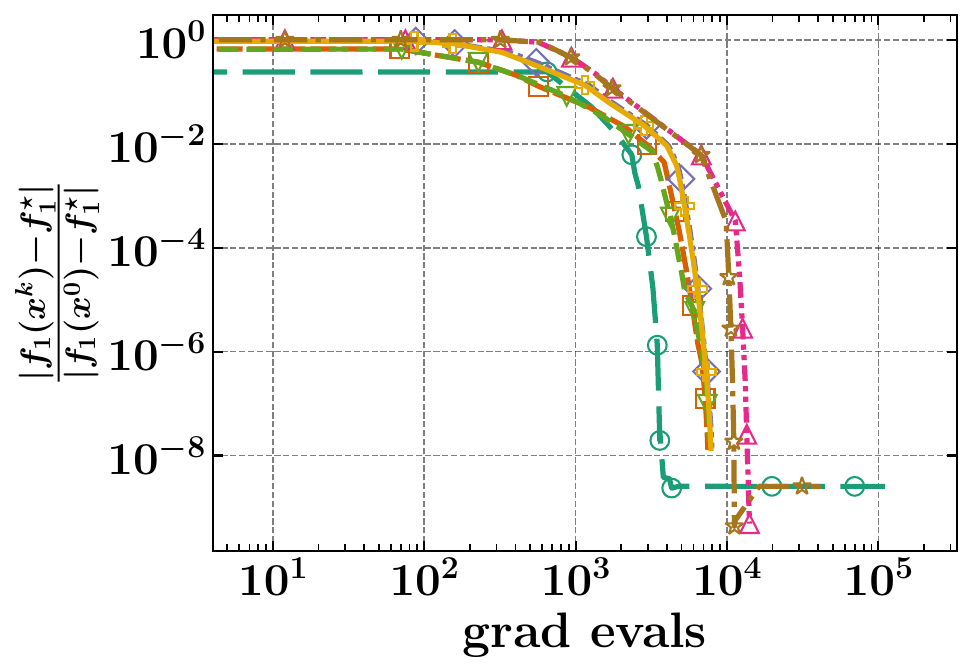}}}
	\caption{Four CUTEst QP instances, one per row, under both inner
	solvers, for $\alpha\in\{1.01,2,8\}$ (and $\xi=10$ for the ALM variant). Each panel plots the
	relative suboptimality against gradient evaluations.}
	\label{fig:ALM-comparison-0-supp}
\end{figure}

The proximal term acts on the conditioning of the subproblems, so its effect is
visible on ill-conditioned instances rather than uniformly. \Cref{fig:dualc1-prox}
isolates it on \texttt{DUALC1}, the worst conditioned instance of the set, with
$\kappa(Q)\approx10^{6}$, under \texttt{PANOC} and at a common $\alpha = 4$.
\texttt{ALM} breaks down. \BALM{} descends but then oscillates
without settling, and does not meet the tolerance within the budget. \PBALM{}
descends monotonically and terminates cleanly. This is the sense in
which the proximal term contributes numerically, namely robustness on
ill-conditioned subproblems rather than a uniform gain in cost. On the
better-conditioned instances of \cref{fig:ALM-comparison-0-supp} the two are
essentially indistinguishable.

\subsection*{Effect of $\alpha$ under a fixed initial penalty}

\begin{figure}[h]\centering
  \resizebox*{9cm}{!}{\includegraphics{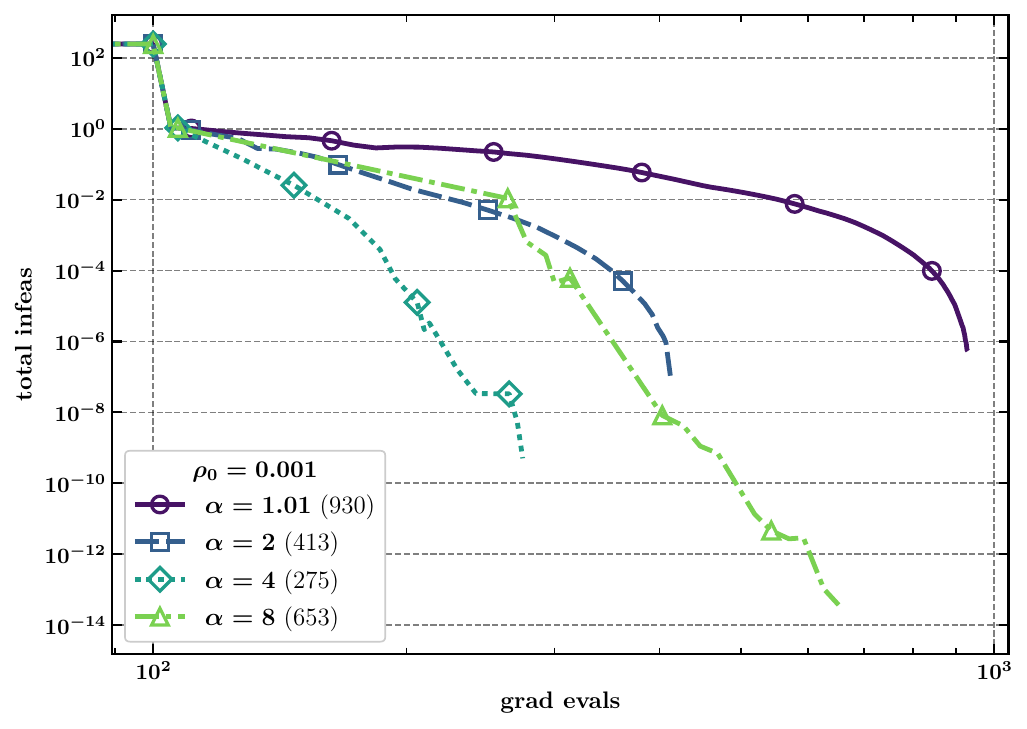}}
  \caption{Sparse nonconvex QCQP under \texttt{PANOC} with $\rho_0$ fixed at
  $10^{-3}$ and $\alpha$ swept, \PBALM{} throughout; gradient evaluations to
  termination in the legend.}
  \label{fig:aux-fixedrho}
\end{figure}

\Cref{fig:aux-fixedrho} holds $\rho_0$ fixed at $10^{-3}$ and sweeps $\alpha$.
In this configuration the schedule is the only mechanism that can raise the
penalty to a useful magnitude, so a larger exponent is needed and the cost is
lowest at an intermediate $\alpha$ rather than at the smallest one. This behavior with respect to the choice of \(\alpha\) is ultimately dictated by the initialization of the penalty parameters.  

\subsection*{Effect of $\alpha$ on CUTEst instances}
The figures below report four CUTEst instances for a wider range of $\alpha$. 
We compare $\alpha\in\{1.01,2,8\}$ against \texttt{ALM} with $\xi=10$, with all initial parameters set as in \cref{sec:setup}, so that $\rho_0$ follows \eqref{eq:rule3} throughout.
On the better conditioned instances every variant converges and the exponent has little
	effect, whereas on \texttt{DUALC1} the runs separate and several fail to
	settle within the budget under either inner solver.
	
\end{document}